\documentclass[leqno,12pt]{article}

\usepackage{amssymb}
\usepackage{euscript}
\usepackage[dvips]{graphicx}
\usepackage{flafter}
\usepackage{pstricks}
\usepackage{pgf,tikz}
\usetikzlibrary{arrows}
\usepackage{caption}
\usepackage[font={footnotesize}]{caption}  

\usepackage{verbatim}

\usepackage{amsmath}

\usepackage{color}

\usepackage{mathrsfs} 

\usepackage{amsthm}

\begin{document}

\baselineskip=17pt
\setcounter{page}{1}

\renewcommand{\theequation}{\thesection.\arabic{equation}}
\newtheorem{theorem}{Theorem}[section]
\newtheorem{lemma}[theorem]{Lemma}
\newtheorem{definition}[theorem]{Definition}
\newtheorem{proposition}[theorem]{Proposition}
\newtheorem{corollary}[theorem]{Corollary}
\newtheorem{fact}[theorem]{Fact}
\newtheorem{problem}[theorem]{Problem}
\newtheorem{conjecture}[theorem]{Conjecture}
\newtheorem{claim}[theorem]{Claim}

\theoremstyle{definition} 
\newtheorem{remark}[theorem]{Remark}
\newtheorem{example}[theorem]{Example}

\newcommand{\eqnsection}{
\renewcommand{\theequation}{\thesection.\arabic{equation}}
    \makeatletter
    \csname  @addtoreset\endcsname{equation}{section}
    \makeatother}
\eqnsection


\def\r{{\mathbb R}}
\def\e{{\mathbb E}}
\def\p{{\mathbb P}}
\def\P{{\bf P}}
\def\E{{\bf E}}
\def\Q{{\bf Q}}
\def\bro{\mathtt{bro}}
\def\z{{\mathbb Z}}
\def\N{{\mathbb N}}
\def\T{{\mathbb T}}
\def\G{G}

\def\ee{\mathrm{e}}
\def\d{\, \mathrm{d}}



\vglue50pt

\centerline{\large\bf The free energy in the Derrida--Retaux recursive model}


\bigskip
\bigskip

\centerline{by}

\medskip

\centerline{Yueyun Hu\footnote{\scriptsize LAGA, Universit\'e Paris XIII, 99 avenue Jean-Baptiste Cl\'ement, F-93430 Villetaneuse, France, {\tt yueyun@math.univ-paris13.fr}} and Zhan Shi\footnote{\scriptsize LPMA, Universit\'e Pierre et Marie Curie, 4 place Jussieu, F-75252 Paris Cedex 05, France, {\tt zhan.shi@upmc.fr}}}

\medskip

\centerline{\it Universit\'e Paris XIII \& Universit\'e Paris VI}

\bigskip
\bigskip
\bigskip

{\leftskip=1truecm \rightskip=1truecm \baselineskip=15pt \small

\noindent{\slshape\bfseries Summary.} We are interested in a simple max-type recursive model studied by Derrida and Retaux~\cite{derrida-retaux} in the context of a physics problem, and find a wide range for the exponent in the free energy in the nearly supercritical regime.

\bigskip

\noindent{\slshape\bfseries Keywords.} Max-type recursive model, free energy.

\bigskip

\noindent{\slshape\bfseries 2010 Mathematics Subject Classification.} 60J80, 82B44.

} 

\bigskip
\bigskip

\section{Introduction}
\label{s:intro}

\subsection{The Derrida--Retaux model}

We are interested in a max-type recursive model investigated in 2014 by Derrida and Retaux \cite{derrida-retaux}, as a toy version of the hierarchical pinning model; see Section \ref{subs:pinning}. The model can be defined, up to a simple change of variables, as follows: for all $n\ge 1$,
\begin{equation}
    X_{n+1}
    \; {\buildrel \mathrm{law} \over =} \;
    (X_n + \widetilde{X}_n -1)^+,
    \label{recurrence_relation}
\end{equation}

\noindent where $\widetilde{X}_n$ denotes an independent copy of $X_n$, and ``${\buildrel \mathrm{law} \over =}$" stands for identity in distribution. We assume that $X_0$ is a non-negative random variable.
    
Since $(X_n + \widetilde{X}_n -1)^+ \le X_n + \widetilde{X}_n$, we have $\e (X_{n+1}) \le 2 \, \e (X_n)$, which implies the existence of the free energy
\begin{equation}
    F_\infty
    :=
    \lim_{n\to \infty} \downarrow \, \frac{\e (X_n)}{2^n} \, .
    \label{F}
\end{equation}

\noindent An immediate question is how to separate the two regimes $F_\infty >0$ and $F_\infty =0$.

\medskip

\begin{example}
\label{ex:bernoulli}

  Assume $\p (X_0 =2) =p$ and $\p (X_0=0) = 1-p$, where $p\in [0, \, 1]$ is a parameter. There exists $p_c \in (0, \, 1)$ such that $F_\infty>0$ if $p>p_c$, and that $F_\infty =0$ if $p<p_c$. 
  
  The value of $p_c$ is known to be $\frac15$ (Collet et al.\ \cite{collet-eckmann-glaser-martin}).\qed

\end{example}

More generally, we write $P_X$ for the law of an arbitrary random variable $X$, and assume from now on
$$
P_{X_0}
=
(1-p) \, \delta_0 + p\, P_{Y_0},
$$

\noindent where $\delta_0$ is the Dirac measure at the origin, $Y_0$ is a random variable taking values in $(0, \, \infty)$, and $p\in [0, \, 1]$ a parameter. In Example \ref{ex:bernoulli}, we have $Y_0=2$.

We often write $F_\infty(p)$ instead of $F_\infty$ in order to make appear the dependence of the free energy in terms of the parameter $p$. Clearly $p\mapsto F_\infty(p)$ is non-decreasing. So there exists a critical parameter $p_c \in [0, \, 1]$ such that
$$
F_\infty(p) > 0 \hbox{ \rm if } p>p_c, 
\qquad
F_\infty(p) = 0 \hbox{ \rm if } p<p_c\, .
$$

\noindent [The extreme cases: $p_c=0$ means $F_\infty(p)>0$ for all $p>0$, whereas $p_c=1$ means $F_\infty(p) =0$ for all $p<1$.] 

We can draw the first $n$ generations of the rooted binary tree leading to the random variable $X_n$; in this sense, $F_\infty(p)$ can be viewed as a kind of percolation function on the binary tree: when $F_\infty>0$, we say there is percolation, whereas if $F_\infty = 0$, we say there is no percolation. {F}rom this point of view, two questions are fundamental: (1) What is the critical value $p_c$? (2) What is the behaviour of the free energy $F_\infty(p)$ when $p$ is in the neighbourhood of $p_c$?

Concerning the first question, the value of $p_c$ can be determined if the random variable $Y_0$ is integer-valued.

\bigskip

\noindent {\bf Theorem A (Collet et al.\ \cite{collet-eckmann-glaser-martin}).} Assume $Y_0$ takes values in $\{ 1, \, 2, \ldots\}$. Then
$$
p_c
=
\frac{1}{\e [(Y_0-1)2^{Y_0}] +1} \, .
$$

\bigskip

Theorem A is proved in \cite{collet-eckmann-glaser-martin} assuming $\e (Y_0\, 2^{Y_0})<\infty$. It is easily seen that it still holds in the case $\e (Y_0\, 2^{Y_0}) = \infty$: Indeed, for $Z_0:= \min\{ Y_0, \, k\}$ in the place of $Y_0$, the corresponding critical value for $p$ is $\frac{1}{\e[(Z_0-1)2^{Z_0}] +1}$, which can be made as close to $0$ as possible by choosing $k$ sufficiently large (by the monotone convergence theorem), so $p_c=0$.

When $Y_0$ is not integer-valued, Theorem A is not valid any more. The value of $p_c$ is unknown (see Section \ref{s:final} for some open problems). However, it is possible to characterise the positivity of $p_c$.

\medskip

\begin{proposition}
\label{p:p_c>0}

 We have $p_c>0$ if and only if $\e(Y_0 \, 2^{Y_0}) <\infty$.

\end{proposition}

\medskip

\noindent {\it Proof.} (1) We first assume that $Y_0$ takes values in $\{ 0, \, 1, \, 2,\ldots\}$.

By Theorem A,
$$
p_c
=
\frac{1}{(\e[(Y_0-1)2^{Y_0}])^+ +1} \, ,
$$

\noindent where $(\e[(Y_0-1)2^{Y_0}])^+$ is the positive part of $\e[(Y_0-1)2^{Y_0}]$. This means that if $\e (Y_0\, 2^{Y_0})<\infty$, then $p_c = \frac{1}{(\e[(Y_0-1)2^{Y_0}])^+ +1}$ (which, in particular, is positive), whereas if $\e (Y_0\, 2^{Y_0})=\infty$, then $p_c =0$.

(2) We now remove the assumption that $Y_0$ is integer-valued. We write
$$
\lfloor Y_0 \rfloor 
\le
Y_0
\le 
\lceil Y_0 \rceil \, .
$$ 

\noindent For both $\lfloor Y_0 \rfloor$ and $\lceil Y_0 \rceil$, we apply the positivity criterion proved in the first step. Since the three conditions $\e(\lfloor Y_0 \rfloor \, 2^{\lfloor Y_0 \rfloor}) <\infty$, $\e(\lceil Y_0 \rceil \, 2^{\lceil Y_0 \rceil}) <\infty$ and $\e(Y_0 \, 2^{Y_0}) <\infty$ are equivalent, the desired result follows.\qed

\bigskip

Proposition \ref{p:p_c>0} tells us that the positivity of $p_c$ does not depend on the exact distribution of $Y_0$, but only on its tail behaviour. 

We now turn our attention to the second question. For the standard Bernoulli bond percolation problem, the percolation function (i.e., the probability that the origin belongs to the unique infinite cluster) is continuous, but not differentiable, at $p=p_c$. For our model, the situation is believed to be very different; in fact, it is predicted (\cite{derrida-retaux}) that the free energy is smooth at $p=p_c$ and that all the derivatives at $p_c$ vanish: 

\medskip

\begin{conjecture}
\label{conj:bernard}

 {\bf (Derrida and Retaux~\cite{derrida-retaux}).} Assume $p_c>0$. There exists a constant $K\in (0, \, \infty)$ such that
 $$
 F_\infty(p)
 =
 \exp\Big( - \frac{K+o(1)}{(p-p_c)^{1/2}} \Big) ,
 \qquad
 p\downarrow p_c \, .
 $$

\end{conjecture}

\medskip

By Proposition \ref{p:p_c>0}, the assumption $p_c>0$ in Conjecture \ref{conj:bernard} means $\e(Y_0 \, 2^{Y_0}) <\infty$.

We have not been able to prove the conjecture. Our aim is to study the influence, on the behaviour of $F_\infty$ near $p_c$, produced by the tail behaviour of $Y_0$. It turns out that our main result can be applied to a more general family of recursive models, which we define in the following paragraph. 

\subsection{A generalised max-type recursive model}

Let $\nu$ be a random variable taking values in $\{ 1, 2, \ldots\}$, such that $m:= \E(\nu) \in (1, \, \infty)$. [We use $\E$ to denote expectation with respect to the law of $\nu$.] For all $n\ge 1$, let
\begin{equation}
    X_{n+1}
    \; {\buildrel \mathrm{law} \over =} \;
    (X_{n,1} + \cdots + X_{n, \nu} -1)^+,
    \label{recurrence_relation_m}
\end{equation}

\noindent where $X_{n,1}$, $X_{n,2}$, $\ldots$ are independent copies of $X_n$, and are independent of $\nu$. {F}rom probabilistic point of view, this is a natural Galton--Watson-type extension of the model in \eqref{recurrence_relation}, which corresponds to the special case $\nu =2$ a.s. We do not know whether there would be a physical interpretation of the randomness of $\nu$ (including in the related models described in Section \ref{subs:pinning} below), as asked by an anonymous referee.

Let $\theta>0$. Let us consider the following situation: There exist constants $0<c_1 \le c_2 <\infty$ such that for all sufficiently large $x$,
\begin{equation}
    c_1 \, \ee^{-\theta x}
    \le
    \p (Y_0 \ge x)
    \le
    c_2\, \ee^{-\theta x} \, .
    \label{power_law}
\end{equation}

\noindent When $\theta>\log m$, we have $p_c>0$ (see Remark \ref{r:Pc>0}; this is in agreement to Proposition \ref{p:p_c>0} if $\nu$ is deterministic); the behaviour of the system in this case is predicted by Conjecture \ref{conj:bernard}. We are interested in the case $\theta \in (0, \, \log m]$. 


\medskip

\begin{theorem}
\label{t:power_law}

 Assume $\E(t^\nu)<\infty$ for some $t>1$, and $m:= \E(\nu)>1$.
 Let $\theta\in (0, \, \log m)$. 
 Under the assumption $\eqref{power_law}$, we have
 $$
 F_\infty(p)
 =
 p^{\beta+o(1)} \, ,
 \qquad
 p \downarrow 0\, ,
 $$
 where $\beta = \beta(\theta) := \frac{\log m}{(\log m)-\theta}$.

\end{theorem}

\medskip

Theorem \ref{t:power_law}, which is not deep, is included in the paper for the sake of completeness. Its analogue in the non-hierarchical pinning setting was known; see \cite{lerouvillois}.


The study of the case $\theta=\log m$ is the main concern of the paper. It turns out that we are able to say more. Fix $\alpha\in \r$. We assume the existence of constants $0<c_3 \le c_4 <\infty$ such that for all sufficiently large $x$,
\begin{equation}
    c_3 \, x^\alpha \, m^{-x}
    \le
    \p (Y_0 \ge x)
    \le
    c_4\, x^\alpha \, m^{-x} \, .
    \label{hyp}
\end{equation}

\noindent 
The main result of the paper is as follows.

\medskip

\begin{theorem}
\label{t:main}

 Let $\alpha>-2$. 
 Assume $\E(t^\nu)<\infty$ for some $t>1$, and $m:= \E(\nu)>1$.
 Under the assumption $\eqref{hyp}$, we have
 $$
 F_\infty(p)
 =
 \exp\Big( - \frac{1}{p^{\chi+o(1)}} \Big) ,
 \qquad
 p\downarrow p_c =0 \, ,
 $$
 where $\chi = \chi(\alpha) := \frac{1}{\alpha+2}$.

\end{theorem}

\medskip

Compared to the original Derrida--Retaux model, additional technical difficulties may appear when $\nu$ is random. For example, the analogue of the fundamental Theorem A is not known (see Problem \ref{prob:pc_GW}).

The proof of the theorem gives slightly more precision: There exists a constant $c_5>0$ such that for all sufficiently small $p>0$,
$$
F_\infty(p)
\le
\exp\Big( - \frac{c_5}{p^{\chi}} \Big) \, .
$$

We will regularly use the following elementary inequalities:
\begin{equation}
    \frac{\e(X_n)-\frac{1}{m-1}}{m^n}
    \le
    F_\infty
    \le
    \frac{\e(X_n)}{m^n} \, ,
    \qquad
    n\ge 1\, .
    \label{encadrement_F}
\end{equation}

\noindent The second inequality follows from \eqref{F}. For the first inequality, it suffices to note that by definition, $\e(X_{n+1}) \ge m \, \e(X_n)-1$, so $n\mapsto \frac{\e(X_n)-\frac{1}{m-1}}{m^n}$ is non-decreasing, and $F_\infty = \lim_{n\to \infty} \uparrow \frac{\e(X_n)-\frac{1}{m-1}}{m^n}$.

An immediate consequence of \eqref{encadrement_F} is the following dichotomy:

$\qquad\bullet$ either $\e(X_n) > \frac{1}{m-1}$ for some $n\ge 1$, in which case $F_\infty >0$;

$\qquad\bullet$ or $\e(X_n) \le \frac{1}{m-1}$ for all $n\ge 1$, in which case $F_\infty=0$.

\subsection{About the Derrida--Retaux model}
\label{subs:pinning}

The Derrida--Retaux model studied in our paper has appeared in several places in both mathematics and physics literatures.

{\bf (a)} The recursion in \eqref{recurrence_relation} belongs to a family of max-type recursive models analysed in the survey paper of Aldous and Bandyopadhyay~\cite{aldous-bandyopadhyay}.

{\bf (b)} The model in \eqref{recurrence_relation} was investigated by Derrida and Retaux \cite{derrida-retaux} to understand the nature of the pinning/depinning transition on a defect line in presence of strong disorder. The problem of the depinning transition has attracted much attention among mathematicians \cite{A1,A2,B1,D2,H1,Ga1,G1,G2,G3,G4,G5,lerouvillois,T1,T2} and physicists \cite{D1,F1,Ga1,KL,M2,M1,Tang} over the last thirty years. Much progress has been made in understanding the question of the relevance of a weak disorder \cite{G1,H1,G4}, i.e., whether a weak disorder is susceptible of modifying the nature of this depinning transition. For strong disorder or even, for a weak disorder when disorder is relevant, it is known that the transition should always be smooth \cite{G5}, but the precise nature of the transition is still controversial \cite{Tang,derrida-retaux,M2}.

It is expected that a similar phase transition should occur in a simplified version of the problem, when the line is constrained to a hierarchical geometry \cite{B2,D1,G3,L1}. Even in this hierarchical version, the nature of the transition is poorly understood. This is why Derrida and Retaux~\cite{derrida-retaux} came up with a toy model which, they argue, should behave like the hierarchical model. This toy model turns out to be sufficiently complicated that many fundamental questions remain open (we include a final section discussing some of these open problems in Section \ref{s:final}).

{\bf (c)} The model in \eqref{recurrence_relation} has also appeared in Collet et al.\ \cite{collet-eckmann-glaser-martin} in their study of spin glass model. 

{\bf (d)} The recursion in \eqref{recurrence_relation} has led to the so-called parking schema; see Goldschmidt and Przykucki~\cite{goldschmidt-przykucki}.  

\bigskip

The rest of the paper and the proofs of the theorems are as follows. In Section \ref{s:outline}, we present a (heuristic) outline the proof of Theorem \ref{t:main}. Section \ref{s:ub} is devoted to the upper bounds in Theorems \ref{t:power_law} and \ref{t:main}. In Section \ref{s:lb}, which is the heart of the paper, we prove the lower bound in Theorem \ref{t:main}. The lower bound in Theorem \ref{t:power_law} is proved in Section \ref{s:power_law}. Finally, we make some additional discussions and present several open problems in Section \ref{s:final}.

\section{Proof of Theorem \ref{t:main}: an outline}
\label{s:outline}

The upper bound in Theorem \ref{t:main}, proved in details in Section \ref{s:ub}, relies on a simple analysis of the moment generating function. The idea of using the moment generating function goes back to Collet et al.\ \cite{collet-eckmann-glaser-martin} and Derrida and Retaux~\cite{derrida-retaux}.

The proof of the lower bound in Theorem \ref{t:main}, quite involving and based on a multi-scale type of argument, is done in two steps. The first step consists of the following estimate: if the initial distribution $X_0$ satisfies, for $t\ge 1$ (say),
$$
\p (X_0 \ge t) \;\; ``\! \ge\! " \;\; p\, t^\alpha\, m^{-t} ,
$$

\noindent then
\begin{equation}
    \p ( X_n \ge t ) \;\; ``\! \ge\! " \;\; \widetilde{p}  \, t^\alpha \, m^{-t}, 
    \label{e:heuristics}
\end{equation}

\noindent for $n\ge 1$ and $t\ge 1$, where
$$
\widetilde{p} := p^2 \,  n^{2+\alpha}\, .
$$

\noindent [See Lemma \ref{l:iteration} for a rigorous statement; for a heuristic explanation of \eqref{e:heuristics}, see below.] This allows us to use inductively the estimate, to arrive at:
$$
\p ( X_j \ge t ) \;\; ``\! \ge\! " \;\; p^a  j^b\, t^\alpha \, m^{-t},
$$

\noindent for all $j$ satisfying $p^a j^b \le 1$, with $a>0$ and $b>0$ such that $b \approx (2+\alpha) a$. By $\e(X_j) = \int_0^\infty \p ( X_j \ge t ) \d t$, we get
$$
\e(X_j)
\;\; ``\! \ge\! " \;\;
\kappa \, p^a j^b \, ,
$$

\noindent where $\kappa>0$ is a (small) constant. [This is Lemma \ref{l:E(Xn)>}.] As such, $\e(X_j) \ge \kappa$, if $j = p^{-a/b}$. We are almost home. The rest of the argument consists in replacing $\kappa$ by a constant greater than $\frac{1}{m-1}$. This is done in the second step.

Let $n\ge 1$. To see why \eqref{e:heuristics} is true, we use a hierarchical representation of the system $(X_i, \, 0\le i\le n)$, due to Collet et al.~\cite{collet-eckmann-glaser-martin2} and Derrida and Retaux~\cite{derrida-retaux}. We define a family of random variables $(X(u), \, u\in \T^{(n)})$, indexed by a reversed Galton--Watson tree $\T^{(n)}$. Let $\T_0 = \T_0^{(n)}$ denote the initial generation of $\T^{(n)}$. We assume that $X(u)$, for $u\in \T_0$, are i.i.d.\ having the distribution of $X_0$. For any vertex $u\in \T^{(n)}\backslash \T_0$, we set
$$
X(u)
:=
(X(u^{(1)}) + \cdots + X(u^{(\nu_u)})-1)^+ \, ,
$$

\noindent where $u^{(1)}$, $\ldots$, $u^{(\nu_u)}$ are the parents of $u$. Consider (see \eqref{Z}; $\lambda_1$ being an unimportant constant, and can be taken as $\frac13$)
$$
Z = Z_n \approx \# \{ u\in \T_0: \, X(u) \ge \lambda_1 n , \, \mathscr{D}(u)\} ,
$$

\noindent where $\mathscr{D}(u) := \{ \exists v \in \T_0 \backslash \{u\} \, , 
 |u\wedge v| < \lambda_1 n, \, X(v) \ge t+n-X(u) \}$, and the symbol ``$\#$" denotes cardinality. Here, by $|u\wedge v| < \lambda_1 n$, we mean $u$ and $v$ have a common descendant before generation $\lambda_1 n$. We observe that $\p(X_n \ge t) \ge \p(Z \ge 1)$. So by the Cauchy--Schwarz inequality,
$$
\p(X_n \ge t)
\ge
\frac{[\e(Z)]^2}{\e(Z^2)} \, .
$$

\noindent The proof of \eqref{e:heuristics} is done by proving that $\e(Z) \approx p^2 n^{2+\alpha} t^\alpha m^{-t}$ (see \eqref{E(Z)bis}) and that $\e(Z^2) \;\; ``\! \le \! " \;\; \e(Z)$.

In the second step, we take $n \approx (\frac1p)^{1/(2+\alpha)}$, defined rigorously in \eqref{n}. We use once again the hierarchical representation. Let $u \in \T_0$ be such that $X(u) = \max_{v\in \T_0} X(v)$. Since the values of $X$ in the initial generation are i.i.d.\ of the law of $X_0$, it is elementary that $X(u) \ge \ell \;\; ``\! :=\! " \;\; n- \log n$ (see \eqref{ell} for the exact value of $\ell$). Let $k:= \frac{\ell}{4}$. The fact $k< \ell$ allows us to use the following inequality:
$$
X_n
\ge
X(u) - n + \sum_{j=0}^{k-1} X_j^{(n)}
\ge
\ell - n + \sum_{j=0}^{k-1} X_j^{(n)} \, ,
$$

\noindent where, for each $j$, $X_j^{(n)}$ has ``approximately" the same law as $X_j$. [For a rigorous formulation, see \eqref{X>}.] Taking expectation on both sides, and using the fact, proved in the first step, that $\e(X_j^{(n)})$ is greater than a constant (say $\kappa_1$) if $j \ge \frac{k}{2}$ (say), we arrive at:
$$
\e(X_n)
\ge
\ell - n + \frac{k}{2} \, \kappa_1 \, ,
$$

\noindent which is greater than a constant multiple of $n$. By the first inequality in \eqref{encadrement_F}, $F_\infty \ge \frac{\e(X_n)-\frac{1}{m-1}}{m^n}$, which implies the desired lower bound in Theorem \ref{t:main}.\qed


\section{Upper bounds}
\label{s:ub}

Without loss of generality, we assume $X_0$ is integer valued (otherwise, we consider $\lceil X_0 \rceil$). Consider the generating functions
$$
G_n(s)
:=
\e(s^{X_n}),
\qquad
h(s)
:=
\E(s^\nu),
$$

\noindent where $\nu$ is the number of independent copies in the convolution relation \eqref{recurrence_relation_m}: $X_{n+1} \; {\buildrel \mathrm{law} \over =} \; (X_{n,1} + \cdots + X_{n, \nu} -1)^+$.

The latter can be written as
$$
G_{n+1}(s)
=
\frac{h(G_n(s))}{s} 
+
\frac{s-1}{s} \, h(G_n(0)) \, .
$$

\noindent Hence
$$ 
G'_{n+1}(s)
=
\frac{h'(G_n(s))}{s} G'_n(s)
-
\frac{h(G_n(s))-h(G_n(0))}{s^2}.
$$

\noindent We fix an $s\in (1, \, m)$ whose value will be determined later. Write
$$
a_n
=
a_n(s)
:=
G_n(s) -1 \, .
$$

\noindent Since $h(G_n(0)) \le h(1) =1$, we have $h(G_n(s))-h(G_n(0)) \ge h(G_n(s))-h(1) = h(1+a_n) - h(1) \ge h'(1) a_n= m \, a_n \ge s\, a_n$. Hence
$$ 
G'_{n+1}(s)
\le
\frac{h'(G_n(s))}{s}  G'_n(s) - \frac{a_n}{s} 
=
\frac{h'(1+a_n)}{s}  G'_n(s) - \frac{a_n}{s} \, .
$$

\noindent By the assumption of existence of $t>1$ satisfying $\E(t^\nu)<\infty$, there exist $\delta_0>0$ and $c_0>0$ such that $h'(1+a) \le h'(1) + c_0 \, a = m + c_0 \, a$ for $a \in (0, \, \delta_0)$. Hence, if $a_n < \delta_0$, then
$$ 
G'_{n+1}(s)
\le
\frac{m + c_0 \, a_n}{s}  G'_n(s) - \frac{a_n}{s} \, .
$$

Let $N_1 := \inf\{n\ge 0: a_n \ge \delta_0 \, \mbox{ or } G_n'(s) \ge \frac{1}{c_0}\}$ (with $\inf \varnothing := \infty$). As long as $a_0< \delta_0$ and $G_0'(s) < \frac{1}{c_0}$ (which we take for granted from now on), we have
$$
    G_{n+1}'(s)
    \le 
    \frac{m +c_0 a_n }{s}  G'_n(s) - \frac{a_n}{s}
    \le 
    \frac{m}{s} G_n'(s),
    \qquad
    0\le n< N_1.
$$

\noindent Iterating the inequality, we get that
\begin{equation}
    G'_n(s)
    \le 
    (\frac{m}{s})^n \, G'_0(s),
    \qquad
    1\le n\le N_1.
    \label{an}
\end{equation}

We now proceed to the proof of the upper bound in Theorem \ref{t:main}. By assumption \eqref{hyp}, $\p (Y_0 \ge x) \le c_4\, x^\alpha \, m^{-x}$ for all sufficiently large $x$. This implies, by integration by parts, the existence of a constant $c_6>0$ such that for all $s\in (1, \, m)$,
$$
\e(Y_0 s^{Y_0-1})
\le
\frac{c_6}{(\log m - \log s)^{\alpha+2}} \, .
$$

\noindent By definition of $G_0$, this yields
\begin{equation}
    G'_0(s) 
    =
    \e( X_0 s^{X_0- 1})
    = 
    p \, \e(Y_0 s^{Y_0-1})
    \le
    \frac{c_6\, p}{(\log m - \log s)^{\alpha+2}} \, ,
    \label{G0_derivative}
\end{equation}

\noindent and for $n\ge 0$,
$$
a_n
:=
G_n(s) -1 
\le 
\e (s^{X_n} \, {\bf 1}_{\{ X_n \ge 1\}}) 
\le 
s\, G'_n(s)
\le
m \, G'_n(s) \, .
$$

\noindent We choose $s:=m \ee^{-1/N}$, where $N := (c_7 \, p)^{-1/(2+\alpha)}$, with $c_7$ denoting a large constant such that $\ee \, \frac{c_6}{c_7} < \frac{1}{c_0}$ and that $\ee\, m \frac{c_6}{c_7} <\delta_0$. Then \eqref{G0_derivative} ensures that $G_0'(s) \le \frac{c_6}{c_7} < \frac{1}{c_0}$, and $a_0\le m\, G'_0(s) < \delta_0$. [In fact, $a_0$ is {\it much} smaller than $m\, G'_0(s)$.]

By \eqref{an}, for $1\le n\le N_1$ (recalling that $G_0'(s) \le \frac{c_6}{c_7}$)
$$
G'_n(s) 
\le 
\ee^{n/N} \, G_0'(s)
\le
\ee^{n/N} \, \frac{c_6}{c_7} \, ,
$$

\noindent and
$$
a_n
\le
m \, G'_n(s)
\le
\ee^{n/N} \, \frac{m \, c_6}{c_7} \, .
$$

\noindent Since $\ee \, \frac{c_6}{c_7} < \frac{1}{c_0}$ and $\ee\, m \frac{c_6}{c_7} <\delta_0$ by the choice of $c_7$, this yields $G'_n(s) < \frac{1}{c_0}$ and $a_n < \delta_0$ for $1\le n\le \min \{ N_1, \, N\}$. By definition of $N_1$, this implies $N_1 \ge N$; hence $\min \{ N_1, \, N\} = N$, which implies $a_N < \delta_0$. 

By Jensen's inequality, $\e(X_n) \le \frac{\log (1+a_n)}{\log s}$. So $\e(X_N) \le \frac{\log (1+ \delta_0)}{\log s}$. In view of the second inequality in \eqref{encadrement_F}, we get, for all sufficiently small $p$,
$$
F_\infty(p)
\le
\frac{\e(X_N)}{m^N}
\le
\frac{\log (1+ \delta_0)}{m^N\, \log s} \, ,
$$

\noindent proving the upper bound in Theorem \ref{t:main}.

The upper bound in Theorem \ref{t:power_law} is obtained similarly. We choose $s:= \ee^{\theta- \varepsilon}$ with $\varepsilon := (\log \frac1{p})^{-1}$. Then
$$
G'_0(s) 
= 
p \, \e(Y_0 s^{Y_0-1})
\le
c_8 \, \frac{p}{\varepsilon^2} \, ,
$$

\noindent for some constant $c_8>0$ and all sufficiently small $p>0$, and $a_0 \le m \, G_0'(s)$. In particular, $G_0'(s) < \frac{1}{c_0}$ and $a_0 < \delta_0$ for all sufficiently small $p>0$. By \eqref{an}, for $1\le n\le N_1$,
$$
G'_n(s) 
\le 
G'_0(s)\, \Big( \frac{m}{\ee^{\theta-\varepsilon}} \Big)^{\! n} 
\le
c_8 \, \frac{p}{\varepsilon^2} \, \Big( \frac{m}{\ee^{\theta-\varepsilon}} \Big)^{\! n} ,
$$

\noindent and
$$
a_n
\le 
m \, G_n'(s)
\le 
m c_8 \, \frac{p}{\varepsilon^2} \, \Big( \frac{m}{\ee^{\theta-\varepsilon}} \Big)^{\! n} \, .
$$

\noindent Let $N' := \frac{\log (c_9 \, \varepsilon^2/p)}{\log (\frac{m}{\ee^{\theta-\varepsilon}} )}$, where $c_9>0$ is a small constant such that $c_8 \, \frac{p}{\varepsilon^2} \, (\frac{m}{\ee^{\theta-\varepsilon}})^{N'} < \frac{1}{c_0}$ and that $m c_8 \, \frac{p}{\varepsilon^2} \, (\frac{m}{\ee^{\theta-\varepsilon}})^{N'} < \delta_0$. Then $G'_n(s) < \frac{1}{c_0}$ and $a_n < \delta_0$ for $1\le n\le \min \{ N_1, \, N'\}$. This implies $a_{N'} < \delta_0$. Since $N'= (1+o(1)) \frac{\log (1/p)}{(\log m) - \theta}$ (for $p\to 0$), we get that $F_\infty(p) \le \frac{\e(X_{N'})}{m^{N'}}\le \frac{\log (1+a_{N'})}{m^{N'} \log s} \le p^{\frac{\log m}{(\log m)- \theta} +o(1)}$, $p\to 0$, proving the upper bound in Theorem \ref{t:power_law}.\qed

\medskip

\begin{remark}
\label{r:Pc>0}

 Let $a_n = a_n(s) := G_n(s) -1$ as in the proof. Since $h(G_n(0)) \le h(1) =1$, we have 
  $$
       a_{n+1}
       = 
       G_{n+1} (s) -1 \le \frac{h(1+a_n) -1}{s} \, .
  $$
 By assumption on $\nu$, there exist $\delta_0'>0$ and $c_0'>0$ such that $h(1+a) -1 \le ma (1+c_0' a)$ for $a \in (0, \, \delta_0')$. Consequently, if $0<a_0<\delta_0'$ and $m(1+c_0'a_0) < s$, then inductively for all $n\ge 0$, 
 $$
 a_{n+1} \le \frac{ma_n (1+c_0' a_n)}{s} < a_n \, .
 $$
 In other words, the sequence $a_n$, $n\ge 1$, is decreasing.
 
 Assume $\p(Y_0>x) \le c_2 \, \ee^{-\theta x}$ for some $\theta>\log m$ and all sufficiently large $x$. Fix $s\in (m, \, \ee^\theta)$. We have $\e(s^{Y_0})<\infty$. So for sufficiently small $p>0$, we have $0<a_0<\delta_0'$ and $m(1+c_0'a_0) \le s$. By the discussions in the last paragraph, the sequence $a_n$, $n\ge 1$, is decreasing. This yields $\sup_{n\ge 0} \e(X_n)<\infty$; hence $F_\infty(p)=0$. Consequently, $p_c>0$ in this case.
 
 [The discussion here gives a sufficient condition for the positivity of $p_c$, not a characterisation.]\qed
 
\end{remark}

\section{Proof of Theorem \ref{t:main}: lower bound}
\label{s:lb}

Throughout the section, we assume $\E(\nu^3) <\infty$ and $m:= \E (\nu)>1$, which is weaker than the assumption in Theorem \ref{t:main}.

We start with a simple hierarchical representation of the system; the idea of this representation already appeared in Collet et al.~\cite{collet-eckmann-glaser-martin2} and in Derrida and Retaux~\cite{derrida-retaux}.

We define a family of random variables $(X(u), \, u\in \T)$, indexed by an infinite tree $\T$, in the following way. For any vertex $u$ in the genealogical tree $\T$, we use $|u|$ to denote the generation of $u$; so $|u|=0$ if $u$ is in the initial generation. We assume that $X(u)$, for $u\in \T$ with $|u|=0$ (i.e., in the initial generation of the system), are i.i.d.\ having the distribution of $X_0$. For any vertex $u\in \T$ with $|u|\ge 1$, let $u^{(1)}$, $\ldots$, $u^{(\nu_u)}$ denote the $\nu_u$ parents of $u$ in generation $|u|-1$, and set
$$
X(u)
:=
(X(u^{(1)}) + \cdots + X(u^{(\nu_u)})-1)^+ \, .
$$

\noindent We assume that for any $i\ge 0$, $(\nu_u, \, |u|=i)$ are i.i.d.\ having the distribution of $\nu$, and are independent of everything else up to generation $i$. 

Fix an arbitrary vertex $\mathfrak{e}_n$ in the $n$-th generation. 
The set of all vertices, including $\mathfrak{e}_n$ itself, in the first $n$ generations of $\T$ having $\mathfrak{e}_n$ as their (unique) descendant at generation $n$, is denoted by $\T^{(n)}$. Clearly, $\T^{(n)}$ is (the reverse of the first $n$ generations of) a Galton--Watson tree, rooted at $\mathfrak{e}_n$, with reproduction distribution $\nu$. Note that $X(\mathfrak{e}_n)$ has the distribution of $X_n$.

More generally, for $v\in \T^{(n)}$ with $|v|=j \le n$, let $\T(v)$ denote the set of all vertices, including $v$ itself, in the first $j$ generations of $\T$ having $v$ as their (unique) descendant at generation $j$. [So $\T^{(n)} = \T (\mathfrak{e}_n)$.] Let $\T_0(v) := \{ x\in \T (v): \, |x|=0\}$. By an abuse of notation, we write $\T_0 := \T_0(\mathfrak{e}_n)$.

Let $v\in \T^{(n)}$. Let $v_{|v|}=v$, and for $|v| <i \le n$, let $v_i$ be the (unique) child of $v_{i-1}$; in particular, $|v_i|=i$ (for $|v| \le i \le n$), and $v_n=\mathfrak{e}_n$. See Figure \ref{f:fig1}.

For $v\in \T^{(n)} \backslash \{ \mathfrak{e}_n\}$, let $\bro (v)$ denote the set of the brothers of $v$, i.e., the set of vertices, different from $v$, that are in generation $|v|$ and having the same child as $v$. Note that $\bro (v)$ can be possibly empty.\footnote{The term ``brothers" is with respect to $\T^{(n)}$, which is {\it reversed} Galton--Watson.}

\begin{figure}[h]

\definecolor{qqqqff}{rgb}{0.,0.,1.}

\centering

\begin{tikzpicture}
[line cap=round,line join=round,>=triangle 45, scale=0.8, every node/.style={transform shape}, x=3.5cm,y=3.5cm]
\clip(-1.4501479165807203,0.42616655413832194) rectangle (2.6940708500405885,2.654566426214505);
\draw (-0.8,2.)-- (-0.6,1.8);
\draw (0.4,2.)-- (0.4,1.8);
\draw [line width=2.8pt] (1.2,1.8)-- (1.2,2.);
\draw (-0.6,1.8)-- (-0.4,1.6);
\draw [line width=2.8pt] (0.6,1.4)-- (0.6,1.2);
\draw [line width=2.8pt] (0.6,1.2)-- (0.6,1.);
\draw [line width=2.8pt] (0.6,1.)-- (0.6,0.8);
\draw [line width=2.8pt] (0.6,0.8)-- (0.6,0.6);
\draw (0.6,1.4)-- (0.4,1.6);
\draw (0.4,1.6)-- (0.4,1.8);
\draw [line width=2.8pt] (0.6,1.4)-- (1.,1.6);
\draw [line width=2.8pt] (1.,1.6)-- (1.2,1.8);
\draw (0.6118782093902231,0.6479982156119691) node[anchor=north west] {$\mathfrak{e}_n$};
\draw (0.25,1.62) node[anchor=north west] {$v_4$};
\draw (0.6118782093902231,1.4496171741190125) node[anchor=north west] {$\mathfrak{e}_5$};
\draw (0.6118782093902231,1.2479520273247877) node[anchor=north west] {$\mathfrak{e}_6$};
\draw (0.6169198380600787,0.8446217337363382) node[anchor=north west] {$\mathfrak{e}_{n-1}$};
\draw (-0.7796113034899245,1.0109954798415737) node[anchor=north west] {$\mathbb{T}^{(n)}$};
\draw (1.0152085029786717,1.6462406922433817) node[anchor=north west] {$\mathfrak{e}_4$};
\draw (1.2168736497728958,1.8428642103677506) node[anchor=north west] {$\mathfrak{e}_3$};
\draw (1.2118320211030402,2.049570985831831) node[anchor=north west] {$\mathfrak{e}_2$};
\draw (-0.3712393812316203,1.6865737216022265) node[anchor=north west] {$w$};
\draw (-0.7947361894994913,2.589025253506382) node[anchor=north west] {\scriptsize$\mathbb{T}_0(w)$};
\draw (-0.8,2.)-- (-1.,2.2);
\draw (-0.8,2.)-- (-0.8,2.2);
\draw (-1.,2.2)-- (-1.,2.4);
\draw (-0.8,2.2)-- (-0.8,2.4);
\draw (-0.4,2.)-- (-0.4,2.2);
\draw (-0.4,2.2)-- (-0.4,2.4);
\draw (0.4,2.)-- (0.4,2.2);
\draw (0.4,2.2)-- (0.4,2.4);
\draw (0.8,2.2)-- (0.8,2.4);
\draw [line width=2.8pt] (1.2,2.)-- (1.2,2.2);
\draw [line width=2.8pt] (1.2,2.2)-- (1.2,2.4);
\draw (1.6,2.)-- (1.6,2.2);
\draw (1.6,2.2)-- (1.6,2.4);
\draw (-1.,2.2)-- (-1.2,2.4);
\draw (1.6,2.2)-- (1.8,2.4);
\draw [rotate around={0.818233648260494:(-0.7618134220909748,2.4242219212226614)},dash pattern=on 4pt off 4pt] (-0.7618134220909748,2.4242219212226614) ellipse (2.2470363378883613cm and 0.4706445145847832cm);
\draw (-0.4,1.6)-- (-0.4,1.8);
\draw (-0.4,2.)-- (-0.4,1.8);
\draw (-0.4,2.2)-- (-0.4,2.4);
\draw (-0.038491889021150257,2.0) node[anchor=north west] {$v$};
\draw (0.25,1.82) node[anchor=north west] {$v_3$};
\draw (1.2118320211030402,2.457942908090136) node[anchor=north west] {$\mathfrak{e}_0$};
\draw (1.2118320211030402,2.2461945039562004) node[anchor=north west] {$\mathfrak{e}_1$};
\draw (1.6151623146914889,2.447859650750425) node[anchor=north west] {$u$};
\draw (1.6151623146914889,2.2461945039562004) node[anchor=north west] {$u_1$};
\draw (1.6151623146914889,2.049570985831831) node[anchor=north west] {$u_2$};
\draw [dash pattern=on 2pt off 2pt] (-0.4,1.6)-- (0.2,0.8);
\draw (0.2,0.8)-- (0.6,0.6);
\draw (0.6,0.6)-- (1.2,0.8);
\draw [dash pattern=on 2pt off 2pt] (1.2,0.8)-- (2.,1.6);
\draw (2.,1.6)-- (2.2,1.8);
\draw (2.2,2.)-- (2.2,1.8);
\draw (2.2,2.)-- (2.2011375552471444,2.1735837295658262);
\draw (2.2011375552471444,2.1735837295658262)-- (2.2,2.4);
\draw (2.2011375552471444,2.1735837295658262)-- (2.4,2.4);
\draw (-0.4,2.2)-- (-0.6,2.4);
\draw (0.4,1.8)-- (0.,2.);
\draw (0.,2.)-- (0.,2.2);
\draw (0.,2.2)-- (0.2,2.4);
\draw (0.4,1.8)-- (0.6,2.);
\draw (0.6,2.)-- (0.6,2.2);
\draw (0.6,2.2)-- (0.6,2.4);
\draw (0.8,2.2)-- (1.,2.4);
\draw (0.,2.2)-- (0.,2.4);
\draw [rotate around={-1.8554662127472683:(0.6026133665684695,1.9998408945810688)},dash pattern=on 4pt off 4pt] (0.6026133665684695,1.9998408945810688) ellipse (1.3229849378156604cm and 0.27614013753876726cm);
\draw (0.6,2.)-- (0.8,2.2);
\draw (0.65,2.08) node[anchor=north west] {\scriptsize$\mathtt{bro}(v)$};
\draw (1.2,1.8)-- (1.6,2.);
\draw (1.6,2.2)-- (2.,2.4);
\begin{scriptsize}
\draw [fill=qqqqff] (-0.8,2.) circle (2.5pt);
\draw [fill=qqqqff] (-0.4,2.) circle (2.5pt);
\draw [fill=qqqqff] (0.4,2.) circle (2.5pt);
\draw [fill=qqqqff] (1.2,2.) circle (2.5pt);
\draw [fill=qqqqff] (1.6,2.) circle (2.5pt);
\draw [fill=qqqqff] (-0.6,1.8) circle (2.5pt);
\draw [fill=qqqqff] (0.4,1.8) circle (2.5pt);
\draw [fill=qqqqff] (1.2,1.8) circle (2.5pt);
\draw [fill=qqqqff] (-0.4,1.6) circle (2.5pt);
\draw [fill=qqqqff] (0.6,0.6) circle (2.5pt);
\draw [fill=qqqqff] (0.6,1.4) circle (2.5pt);
\draw [fill=qqqqff] (0.6,1.2) circle (2.5pt);
\draw [fill=qqqqff] (0.6,1.) circle (2.5pt);
\draw [fill=qqqqff] (0.6,0.8) circle (2.5pt);
\draw [fill=qqqqff] (0.4,1.6) circle (2.5pt);
\draw [fill=qqqqff] (1.,1.6) circle (2.5pt);
\draw [fill=qqqqff] (-1.,2.2) circle (2.5pt);
\draw [fill=qqqqff] (-0.8,2.2) circle (2.5pt);
\draw [fill=qqqqff] (-1.,2.4) circle (2.5pt);
\draw [fill=qqqqff] (-0.8,2.4) circle (2.5pt);
\draw [fill=qqqqff] (-0.6,2.4) circle (2.5pt);
\draw [fill=qqqqff] (-0.4,2.2) circle (2.5pt);
\draw [fill=qqqqff] (-0.4,2.4) circle (2.5pt);
\draw [fill=qqqqff] (0.2,2.4) circle (2.5pt);
\draw [fill=qqqqff] (0.4,2.2) circle (2.5pt);
\draw [fill=qqqqff] (0.4,2.4) circle (2.5pt);
\draw [fill=qqqqff] (0.8,2.2) circle (2.5pt);
\draw [fill=qqqqff] (0.8,2.4) circle (2.5pt);
\draw [fill=qqqqff] (1.,2.4) circle (2.5pt);
\draw [fill=qqqqff] (1.2,2.2) circle (2.5pt);
\draw [fill=qqqqff] (1.2,2.4) circle (2.5pt);
\draw [fill=qqqqff] (1.6,2.2) circle (2.5pt);
\draw [fill=qqqqff] (1.6,2.4) circle (2.5pt);
\draw [fill=qqqqff] (1.8,2.4) circle (2.5pt);
\draw [fill=qqqqff] (2.,2.4) circle (2.5pt);
\draw [fill=qqqqff] (2.2,2.4) circle (2.5pt);
\draw [fill=qqqqff] (-1.2,2.4) circle (2.5pt);
\draw [fill=qqqqff] (2.4,2.4) circle (2.5pt);
\draw [fill=qqqqff] (0.6,2.4) circle (2.5pt);
\draw [fill=qqqqff] (-0.4,1.8) circle (2.5pt);
\draw [fill=qqqqff] (0.2,0.8) circle (2.5pt);
\draw [fill=qqqqff] (2.2,1.8) circle (2.5pt);
\draw [fill=qqqqff] (1.2,0.8) circle (2.5pt);
\draw [fill=qqqqff] (2.,1.6) circle (2.5pt);
\draw [fill=qqqqff] (2.2,2.) circle (2.5pt);
\draw [fill=qqqqff] (2.2011375552471444,2.1735837295658262) circle (2.5pt);
\draw [fill=qqqqff] (0.,2.) circle (2.5pt);
\draw [fill=qqqqff] (0.,2.2) circle (2.5pt);
\draw [fill=qqqqff] (0.6,2.) circle (2.5pt);
\draw [fill=qqqqff] (0.6,2.2) circle (2.5pt);
\draw [fill=qqqqff] (0.,2.4) circle (2.5pt);
\end{scriptsize}
\end{tikzpicture}
\caption{Spine $(\mathfrak{e}_i)_{0\le i \le n}$;  in the above picture, $v_2:=v$ and $v_i=\mathfrak{e}_i$ for all $5\le i \le n.$}
\label{f:fig1}

\end{figure}

In the rest of the paper, we use $P_\omega$ the conditional probability given $\T$, and its corresponding expectation $E_\omega$. The law of $\T$ is denoted by $\P$, the corresponding expectation $\E$. We write $\p (\, \cdot \,) := \E [ P_\omega (\, \cdot \,)]$, with corresponding expectation $\e$.

We now describe the law of the {\it size-biased Galton--Watson tree}.
Let $\Q$ be the probability measure defined on $\sigma(\T^{(n)})$, the sigma-field generated by $\T^{(n)}$, by $\Q := \frac{\# \T_0}{m^n} \bullet \P$. Under $\Q$, $\T^{(n)}$ represents (the first $n$ generations of) a so-called size-biased Galton--Watson tree. There is a simple way to describe the law of the size-biased Galton--Watson tree. Let $\mathfrak{e}_0 = \mathfrak{e}_0^{(n)}$ be a random variable taking values in $\T_0$ (which is not measurable with respect to $\sigma(\T_0)$, the sigma-field generated by $\T_0$) whose under $\Q$, given $\sigma(\T_0)$, is uniformly distributed on $\T_0$: $\Q( \mathfrak{e}_0 = u \, | \, \T_0) = \frac{1}{\# \T_0}$ for any $u\in \T_0$. Let $\mathfrak{e}_i= \mathfrak{e}_i^{(n)}$ be the unique descendant at generation $i$ of $\mathfrak{e}_0$, for all $0\le i\le n$. The collection $(\mathfrak{e}_i, \, 0\le i\le n)$ is referred to as the spine. The {\it spinal decomposition theorem} says that under $\Q$, $\bro (\mathfrak{e}_i)$, for $0\le i\le n-1$, are i.i.d., and conditionally on $\mathscr{G}_n := \sigma (\bro (\mathfrak{e}_i), \, 0\le i\le n-1)$, $(\T_0(v), \, v\in \bro (\mathfrak{e}_i))_{0\le i\le n-1}$ are independent; furthermore, for each $0\le i\le n-1$, conditionally on $\mathscr{G}_n$ and under $\Q$, $\T_0(v)$, for $v\in \bro (\mathfrak{e}_i)$ are i.i.d.\ having the law of $\T_0 (u_i)$ under $\P$ (for any $u\in \T_0$ which is $\sigma(\T_0)$-measurable). For more details, see Lyons, Pemantle and Peres~\cite{lyons-pemantle-peres}, or Lyons and Peres (\cite{lyons-peres}, Chap.~12), Shi (\cite{stf}, Chap.~2) for formalism for tree-valued random variables. 

A useful consequence of the spinal decomposition theorem is the many-to-one formula: For any measurable function $g$,
\begin{equation}
    \E \Big( \sum_{u\in \T_0} g(\T(u_1), \ldots, \T(u_n)) \Big)
    =
    m^n \, \E_\Q [ g(\T(\mathfrak{e}_1), \ldots, \T(\mathfrak{e}_n)) ] ,
    \label{many-to-one}
\end{equation}

\noindent where $(\mathfrak{e}_i, \, 0\le i\le n)$ is the spine.

Here is another consequence of the spinal decomposition theorem. Let $1\le i\le n$. We have $\E_\Q [\# \T_0(\mathfrak{e}_i)] = \E_\Q [\# \T_0(\mathfrak{e}_{i-1})] + c_{10} \, m^{i-1}$, where $c_{10} := \E_\Q [ \# \bro (\mathfrak{e}_{i-1}) ] = \E_\Q [ \# \bro (\mathfrak{e}_0) ] = \frac1m \, \sum_{k=1}^\infty k(k-1) \P (\nu =k) <\infty$. This yields
\begin{equation}
    \E_\Q [\# \T_0(\mathfrak{e}_i)]
    =
    c_{10} \, \sum_{j=0}^{i-1} m^j +1
    \le
    c_{11} \, m^i \, ,
    \label{E(T(ek))}
\end{equation}

\noindent for some constant $c_{11}>0$ and all $0\le i\le n$. Similarly, the assumption $\E(\nu^3) <\infty$ ensuring $\E_\Q [ (\# \bro (\mathfrak{e}_0))^2 ] = \frac1m \, \sum_{k=1}^\infty k(k-1)^2 \P (\nu =k) <\infty$, we have
\begin{equation}
    \E_\Q \{ [\# \T_0(\mathfrak{e}_i)]^2 \}
    \le
    c_{12}\, m^{2i} \, ,
    \label{E(T(ek)2)}
\end{equation}

\noindent for some constant $c_{12}>0$ and all $0\le i\le n$.

We now turn to the proof of the lower bound in Theorem \ref{t:main}, which is done in two steps. The first step, summarised in Lemma \ref{l:iteration} below, is a probability estimate that allows for iteration. The second step says that along the spine, $X_n$ will reach sufficiently high expected values.

\subsection{First step: Inductive probability estimate}

The first step gives a useful inductive probability estimate. In order to make the induction possible, we assume something more general than the assumption \eqref{hyp} in Theorem \ref{t:main}. 

\medskip

\begin{lemma} 
\label{l:iteration}

 Assume $\E(\nu^3) <\infty$ and $m:= \E (\nu)>1$.
 Let $\alpha \in \r$. Let $c_{13}>0$, $c_{14}>0$ and $c_{15}>0$.
 There exists $c>0$ such that for $0< p< 1$ with $p \, n^{1+\alpha} \le c_{13}$, if the initial distribution of $X_0$ is such that for some $1\le \gamma \le c_{14}\, n$,
 \begin{equation}
     \p (X_0 \ge t) \ge c_{15} \, p\, (t+\gamma)^\alpha\, m^{-t}, 
     \qquad 
     \forall\,  t>0, 
     \label{gamma}
 \end{equation}

 \noindent then
 $$ 
 \p ( X_n \ge t ) \ge c \, p^2 \,  n^{2+\alpha} \, (t+n)^\alpha \, m^{-t}, 
 \qquad 
 \forall \,  t>0.
 $$ 

\end{lemma}

\medskip

\noindent {\it Proof.} Without loss of generality, we assume $X_0$ is integer valued such that
\begin{equation}
    \p(X_0 = t) = c_{15} \, p\, (t+\gamma)^\alpha\, m^{-t}, \qquad \forall\,  t= 1, \, 2, \ldots
    \label{P(X_0=t)}
\end{equation}

\noindent [In fact, the distribution in \eqref{P(X_0=t)} is stochastically smaller than or equal to a distribution satisfying \eqref{gamma}, with a possibly different value of the constant $c_{15}$.]

For $v\in \T^{(n)}$ with $|v|=j \le n$, let
$$
M(v)
:=
\max_{r \in \bro (v)} \max_{w \in \T_0(r)} (X(w)-j)^+,
$$

\noindent where $r \in \bro (v)$ means, as before, that $r$ is a brother of $v$, and $X(w)$ is the random variable assigned to the vertex $w$ on the initial generation. Let $b>0$. Let $0<\lambda_1<\lambda_2<1$ be any fixed constants.\footnote{The values of $\lambda_1$ and $\lambda_2$ play no significant role in the proof; so we can take, for example, $\lambda_1= \frac13$ and $\lambda_2 = \frac23$.} We consider the integer-valued random variable
\begin{equation}
    Z
    :=
    \sum_{u\in \T_0} {\bf 1}_{A_u} ,
    \label{Z}
\end{equation}

\noindent where, for any $u\in \T_0$,
$$
A_u
:=
\Big\{ X(u) \in [\lambda_1 n, \, \lambda_2 n], \; \max_{1\le j\le \lambda_1 n} M(u_j) \ge b +n- X(u) \Big\} .
$$

\noindent Clearly,
\begin{equation}
    \p(Z\ge 1)  \le \p(X_n \ge b).
    \label{lb}
\end{equation}

\noindent Throughout the proof, we write $x \lesssim y$ or $y\gtrsim x$ if $x\le c y$ for some constant $c\in (0, \, \infty)$ that does not depend on $(n,\, p, \, b)$, and $x\asymp y$ if both relations $x \lesssim y$ and $y\lesssim x$ hold.

For $x \ge (1-\lambda_2)n$ and $1\le j\le n$, we have, by \eqref{P(X_0=t)},
$$
\p (X_0 \ge x+j)
\asymp
p(x+j+\gamma)^\alpha m^{-x-j}
\asymp
px^\alpha m^{-x-j}.
$$

\noindent [So the parameter $\gamma$ figuring in the condition \eqref{gamma} disappears because $x+\gamma \asymp x$ if $x \ge (1-\lambda_2)n$.] Note that $M(u_j)$, for $1\le j\le n$, are independent under $P_\omega$. We have, for $u\in \T_0$, $x \ge (1-\lambda_2)n$, and integers $1\le n_1 \le n_2 \le n$,
\begin{eqnarray*}
    P_\omega \Big( \max_{n_1\le j\le n_2} M(u_j) \ge x \Big)
 &=& 1 - \prod_{j=n_1}^{n_2} P_\omega (M(u_j) < x)
    \\
 &=& 1 - \prod_{j=n_1}^{n_2} \prod_{r \in \bro (u_j)} \prod_{w \in \T_0(r)} [1-\p (X_0 \ge x+j)]
    \\
 &\asymp& \Big[ p x^\alpha \sum_{j=n_1}^{n_2} m^{-x-j} \Lambda (u_j) \Big] \wedge 1,
\end{eqnarray*}

\noindent uniformly in $1\le n_1 \le n_2 \le n$. Here, $a\wedge b := \min\{ a, \, b\}$ for real numbers, and for all $v\in \T^{(n)}$,
$$
\Lambda (v)
:=
\sum_{r\in \bro (v)} \#\T_0 (r) \, ,
$$

\noindent with $\#\T_0 (r)$ denoting the cardinality of $\T_0 (r)$. For future use, we observe that
\begin{equation}
    P_\omega \Big( \max_{n_1\le j\le n_2} M(u_j) \ge x \Big)
    \lesssim
    p x^\alpha \sum_{j=n_1}^{n_2} m^{-x-j} \Lambda (u_j),
    \label{M}
\end{equation}

\noindent uniformly in $1\le n_1 \le n_2 \le n$. Taking $n_1 :=1$ and $n_2 := \lambda_1 n$, and by independence of $X(u)$ and $\max_{n_1\le j\le n_2} M(u_j)$ under $P_\omega$, we arrive at:\footnote{For notational simplification, we treat $\lambda_1 n$ and $\lambda_2 n$ as if they were integers.}
\begin{equation}
    P_\omega (A_u)
    \asymp
    \sum_{\ell=\lambda_1 n}^{\lambda_2 n} \, p n^\alpha m^{-\ell} \, \Big\{ \Big[ p (b+n)^\alpha \sum_{j=1}^{\lambda_1 n} m^{-(b+n-\ell)-j} \Lambda (u_j) \Big] \wedge 1 \Big\} .
    \label{P(Au)=}
\end{equation}

\noindent [We have used $\ell+\gamma \asymp n$ and $b+n-\ell \asymp b+n$.] For future use, we see that by removing ``$\wedge 1$" on the right-hand side,
\begin{equation}
    P_\omega (A_u)
    \lesssim
    p^2 n^{1+\alpha} \, (b+n)^\alpha \sum_{j=1}^{\lambda_1 n} m^{-(b+n)-j} \Lambda (u_j) .
    \label{P(Au)<}
\end{equation}

We now estimate $\e(Z)$ and $\e(Z^2)$.

We first look at the expectation of $Z$ under $P_\omega$: By \eqref{P(Au)=},
$$
E_\omega(Z)
\asymp
\sum_{u\in \T_0} \sum_{\ell=\lambda_1 n}^{\lambda_2 n} \, p \, n^\alpha m^{-\ell} \, \Big\{ \Big[ p (b+n)^\alpha \sum_{j=1}^{\lambda_1 n} m^{-b-n+\ell-j} \Lambda (u_j) \Big] \wedge 1 \Big\} .
$$

\noindent We take expectation on both sides with respect to $\P$, the law of $\T$. By the many-to-one formula \eqref{many-to-one},
\begin{equation}
    \e(Z) 
    \asymp
    m^n \sum_{\ell=\lambda_1 n}^{\lambda_2 n} \, p \, n^\alpha m^{-\ell} \, \E_\Q (\eta \wedge 1) ,
    \label{E(Z)}
\end{equation}

\noindent where
$$
    \eta
    =
    \eta(n, \, \ell)
    := p (b+n)^\alpha \sum_{j=1}^{\lambda_1 n} m^{-b-n+\ell} \, m^{-j} \Lambda (\mathfrak{e}_j) \, . 
$$

\noindent By the spinal decomposition theorem, under $\Q$, $m^{-j} \Lambda (\mathfrak{e}_j) 
$, $1\le j\le n$ are independent, and for each $j$, $\# \bro (\mathfrak{e}_j)$ has the same law as $\#\bro (\mathfrak{e}_1)$, whereas conditionally on $\# \bro (\mathfrak{e}_j)$, $\Lambda (\mathfrak{e}_j)$ is distributed as the sum of $\# \bro (\mathfrak{e}_j)$ independent copies of $\T_0 (u_j)$ under $\P$ (for any $u\in \T_0$), the latter being the number of individuals in the $j$-th generation of a Galton--Watson process with reproduction law $\nu$ (starting with $1$ individual). Accordingly,
$$
\E_\Q [m^{-j} \Lambda (\mathfrak{e}_j)\, | \, \#\bro (\mathfrak{e}_j)] 
=
\#\bro (\mathfrak{e}_j) \, ,
$$

\noindent so that
$$
\E_\Q [m^{-j} \Lambda (\mathfrak{e}_j)]
=
\E_\Q[\#\bro (\mathfrak{e}_j)]
=
\E_\Q[\#\bro (\mathfrak{e}_1)]
<\infty \, .
$$

\noindent Hence
\begin{equation}
    \E_\Q (\eta)
    \asymp
    p (b+n)^\alpha \sum_{j=1}^{\lambda_1 n} m^{-b-n+\ell} 
    \asymp 
    p (b+n)^\alpha n\, m^{-b-n+\ell}\, .
    \label{E(Z)0}
\end{equation}

We now estimate $\E_\Q(\eta\wedge 1)$. Consider a Galton--Watson process with reproduction law $\nu$ (starting with $1$ individual) under $\P$. For each $j\ge 0$, let $m^j W_j$ denote the number of individuals in the $j$-th generation. By Athreya and Ney~\cite{athreya-ney} (p.~9, Theorem 2), as long as $\nu$ has a finite second moment, $(W_j, \, j\ge 0)$ is a martingale bounded in $L^2$; in particular, $W_j$ converges in $L^2$, when $j\to \infty$, to a limit denoted by $W$. For any $s>0$, 
$$
\E_\Q ( \ee^{-s m^{-j} \Lambda (\mathfrak{e}_j)}\, | \, \#\bro (\mathfrak{e}_j)) 
=
[ \E( \ee^{-s W_j} ) ]^{\#\bro (\mathfrak{e}_j)} ,
$$

\noindent so that
$$
\E_\Q ( \ee^{-s m^{-j} \Lambda (\mathfrak{e}_j)} ) 
=
\E_\Q \Big\{ [ \E( \ee^{-s W_j} ) ]^{\#\bro (\mathfrak{e}_j)} \Big\}
=
\E_\Q \Big\{ [ \E( \ee^{-s W_j} ) ]^{\#\bro (\mathfrak{e}_0)} \Big\} .
$$

\noindent By conditional Jensen's inequality, $\ee^{-s W_j} = \ee^{-s \E(W \, | \, W_j)} \le \E(\ee^{-s W} \, | \, W_j)$, so $\E( \ee^{-s W_j} ) \le \E(\ee^{-s W})$. Since $\E(W) =1$, we have $1- \E(\ee^{-s W}) \sim s$, $s\to 0$. Hence, there exists a constant $c_{16}>0$ such that for all $s\in (0, \, 1]$ and all $j\ge 0$,
$$
\E_\Q ( \ee^{-s m^{-j} \Lambda (\mathfrak{e}_j)} ) 
\le
\E_\Q \Big\{ [1-c_{16}\, s ]^{\#\bro (\mathfrak{e}_0)} \Big\}
\le
1- c_{17}\, c_{16}\, s ,
$$

\noindent with $c_{17} := \Q(\#\bro (\mathfrak{e}_0) \ge 1) >0$. Consequently, with $c_{18} := c_{17}\, c_{16}$,
\begin{eqnarray*}
    \E_\Q (\eta \wedge 1)
 &\asymp& 1-\E_\Q (\ee^{-\eta})
    \\
 &=& 1- \prod_{j=1}^{\lambda_1 n} \E_\Q \Big( \ee^{-p (b+n)^\alpha \, m^{-b-n+\ell} \, m^{-j} \Lambda (\mathfrak{e}_j)} \Big)
    \\
 &\ge& 1- \prod_{j=1}^{\lambda_1 n} \Big( 1- c_{18} \, p (b+n)^\alpha \, m^{-b-n+\ell} \Big) 
    \\
 &\asymp&p (b+n)^\alpha \, n \, m^{-b-n+\ell} \, .
\end{eqnarray*}

\noindent Going back to \eqref{E(Z)0}, this yields $\E_\Q (\eta \wedge 1) \asymp p (b+n)^\alpha \, n \, m^{-b-n+\ell}$. In view of \eqref{E(Z)}, we obtain:
\begin{equation}
    \e(Z) 
    \asymp
m^n \sum_{\ell=\lambda_1 n}^{\lambda_2 n} \, p \, n^\alpha m^{-\ell} \, p (b+n)^\alpha \, n \, m^{-b-n+\ell}
    \asymp
    p^2 \, n^{2+\alpha} (b+n)^\alpha \, m^{-b} \, ,
    \label{E(Z)bis}
\end{equation}

\noindent uniformly in $b\ge 0$. 

For the second moment of $Z$, we write, by an abuse of notation, 
$$
\T_k := \{ x\in \T^{(n)}: \, |x|=k \},
\qquad
0\le k\le n;
$$

\noindent then
$$
E_\omega(Z^2)
=
E_\omega(Z) + \sum_{k=1}^n \; \sum_{x\in \T_k} \; \sum_{(u, \, v)} P_\omega(A_u \cap A_v) ,
$$

\noindent where $\sum_{(u, \, v)}$ is over the pairs $(u, \, v)\in \T_0 \times \T_0$ with $u_k=v_k=x$ such that $u_{k-1} \not= v_{k-1}$. We take expectation with respect to $\P$ on both sides, while splitting the sum $\sum_{k=1}^n$ into $\sum_{k= \lambda_1 n+1}^n$ and $\sum_{k=1}^{\lambda_1 n}$:
\begin{eqnarray}
    \e(Z^2)
 &=& \e(Z) 
    +
    \sum_{k= \lambda_1 n+1}^n \E \sum_{x\in \T_k} \; \sum_{(u, \, v)} P_\omega(A_u \cap A_v) 
    \nonumber
    \\
 && + 
    \sum_{k=1}^{\lambda_1 n} \E \sum_{x\in \T_k} \; \sum_{(u, \, v)} P_\omega(A_u \cap A_v) .
    \label{E(Z2)}
\end{eqnarray}

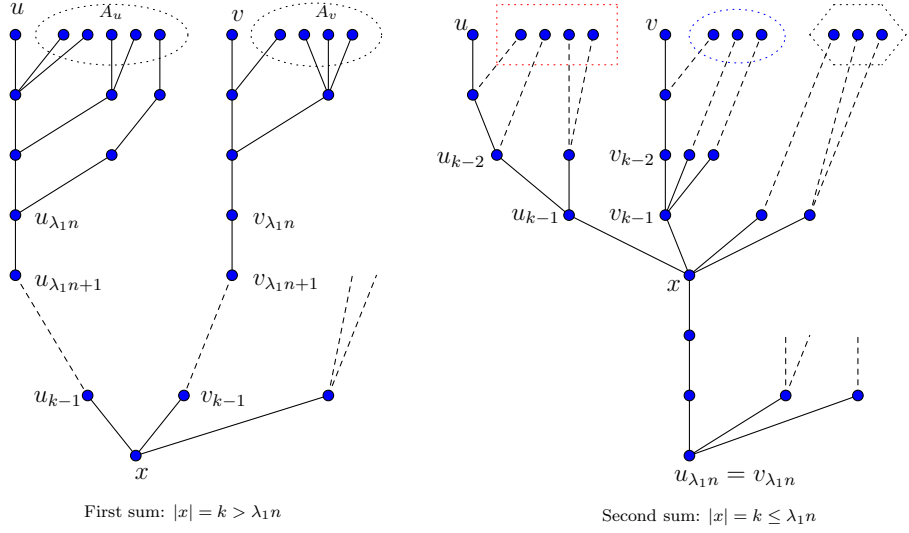
\begin{figure}[h]
\definecolor{ffqqqq}{rgb}{1.,0.,0.}
\definecolor{qqqqff}{rgb}{0.,0.,1.}
\begin{tikzpicture}[line cap=round,line join=round,>=triangle 45, scale=0.8, every node/.style={transform shape}, x=2.0cm,y=2.5cm]
\clip(-1.9012989954675643,-1.5379074448830765) rectangle (6.692498720796524,2.5275878879275417);
\draw (-1.2,2.)-- (-1.2,1.8);
\draw (-1.2,1.8)-- (-1.2,1.6);
\draw (-1.2,1.6)-- (-1.2,1.4);
\draw (-1.2,1.4)-- (-1.2,1.2);
\draw (-1.2,1.2)-- (-1.2,1.);
\draw (-1.2,1.)-- (-1.2,0.8);
\draw (-0.4,1.6)-- (-0.4,2.);
\draw (-0.4,1.6)-- (-0.2,2.);
\draw (-1.2,1.2)-- (-0.4,1.6);
\draw (-1.2,0.8)-- (-0.4,1.2);
\draw (-1.2,1.6)-- (-0.8,2.);
\draw (-1.2,1.6)-- (-0.6,2.);
\draw [rotate around={0.:(-0.4,2.)},dotted] (-0.4,2.) ellipse (1.2661380654259364cm and 0.504829180144908cm);
\draw (-1.2,0.8)-- (-1.2,0.4);
\draw (0.6,0.4)-- (0.6,0.8);
\draw (0.6,0.8)-- (0.6,1.2);
\draw (0.6,1.2)-- (0.6,1.6);
\draw (0.6,1.6)-- (0.6,2.);
\draw (0.6,1.6)-- (1.,2.);
\draw (0.6,1.2)-- (1.4,1.6);
\draw (1.4,1.6)-- (1.2,2.);
\draw (1.4,1.6)-- (1.4,2.);
\draw (1.4,1.6)-- (1.6,2.);
\draw [rotate around={0.:(1.3,2.)},dotted] (1.3,2.) ellipse (1.084819643987329cm and 0.5256449312226121cm);
\draw [dash pattern=on 2pt off 2pt] (-1.2,0.4)-- (-0.6,-0.4);
\draw [dash pattern=on 2pt off 2pt] (0.6,0.4)-- (0.2,-0.4);
\draw (-0.2,-0.8)-- (-0.6,-0.4);
\draw (-0.2,-0.8)-- (0.2,-0.4);
\draw (-0.2,-0.8)-- (1.4,-0.4);
\draw (-0.7002066943947333,-1.063730469718411) node[anchor=north west] {\scriptsize$\mbox{First sum: } |x|=k> \lambda_1 n$};
\draw (-1.12,-0.31519277898227865) node[anchor=north west] {$u_{k-1}$};
\draw (0.26665448947277315,-0.31519277898227865) node[anchor=north west] {$v_{k-1}$};
\draw (-1.1134618778219738,0.4567367145893578) node[anchor=north west] {$u_{\lambda_1 n+1}$};
\draw (-1.1134618778219738,0.8699918980165975) node[anchor=north west] {$u_{\lambda_1 n}$};
\draw (0.7,0.46453398220119246) node[anchor=north west] {$v_{\lambda_1 n+1}$};
\draw (0.7,0.8699918980165975) node[anchor=north west] {$v_{\lambda_1 n}$};
\draw (-0.5754504126053775,2.25) node[anchor=north west] {\scriptsize$A_u$};
\draw (1.2179211381166102,2.25) node[anchor=north west] {\scriptsize$A_v$};
\draw (-1.3161908357296768,2.265702800535011) node[anchor=north west] {$u$};
\draw (0.5239643206633192,2.2189191948640024) node[anchor=north west] {$v$};
\draw (4.4,-0.8)-- (4.4,-0.4);
\draw (4.4,-0.4)-- (4.4,0.);
\draw (4.4,0.)-- (4.4,0.4);
\draw (4.4,0.4)-- (3.4,0.8);
\draw (4.4,0.4)-- (4.2,0.8);
\draw (4.4,0.4)-- (5.,0.8);
\draw (4.4,0.4)-- (5.4,0.8);
\draw (3.4,0.8)-- (2.8,1.2);
\draw (3.4,0.8)-- (3.4,1.2);
\draw (4.2,0.8)-- (4.2,1.2);
\draw (4.2,0.8)-- (4.4,1.2);
\draw (2.8,1.2)-- (2.6,1.6);
\draw (2.6,1.6)-- (2.6,2.);
\draw [dash pattern=on 2pt off 2pt] (2.6,1.6)-- (3.,2.);
\draw [dash pattern=on 2pt off 2pt] (2.8,1.2)-- (3.2,2.);
\draw [dash pattern=on 2pt off 2pt] (3.4,1.2)-- (3.4,2.);
\draw (4.2,1.2)-- (4.2,1.6);
\draw (4.2,1.6)-- (4.2,2.);
\draw (4.2,0.8)-- (4.6,1.2);
\draw [dash pattern=on 2pt off 2pt] (4.4,1.2)-- (4.8,2.);
\draw [dash pattern=on 2pt off 2pt] (4.6,1.2)-- (5.,2.);
\draw [dash pattern=on 2pt off 2pt] (4.2,1.6)-- (4.6,2.);
\draw [dash pattern=on 2pt off 2pt] (5.,0.8)-- (5.6,2.);
\draw [dash pattern=on 2pt off 2pt] (5.4,0.8)-- (5.8,2.);
\draw [dash pattern=on 2pt off 2pt] (5.4,0.8)-- (6.,2.);
\draw (4.14,0.42334491175385365) node[anchor=north west] {$x$};
\draw (-0.27915424335565786,-0.8220151737515349) node[anchor=north west] {$x$};
\draw (2.3719167446681495,2.1643383215811594) node[anchor=north west] {$u$};
\draw (3.9625593374824346,2.1643383215811594) node[anchor=north west] {$v$};
\draw (2.85,0.9011809684639364) node[anchor=north west] {$u_{k-1}$};
\draw (3.66,0.9011809684639364) node[anchor=north west] {$v_{k-1}$};
\draw (2.23,1.2754498138320025) node[anchor=north west] {$u_{k-2}$};
\draw (3.66,1.2910443490556718) node[anchor=north west] {$v_{k-2}$};
\draw (4.2,-0.8376097089752044) node[anchor=north west] {$u_{\lambda_1 n}=v_{\lambda_1 n}$};
\draw (3.6038850273380367,-1.0949195401657499) node[anchor=north west] {\scriptsize$\mbox{Second sum: } |x|=k \le \lambda_1 n$};
\draw (4.4,-0.8)-- (5.8,-0.4);
\draw (4.4,-0.8)-- (5.2,-0.4);
\draw [dash pattern=on 2pt off 2pt] (5.2,-0.4)-- (5.2,0.);
\draw [dash pattern=on 2pt off 2pt] (5.2,-0.4)-- (5.4,0.);
\draw [dash pattern=on 2pt off 2pt] (5.8,-0.4)-- (5.8,0.);
\draw [dash pattern=on 2pt off 2pt] (3.4,1.2)-- (3.6,2.);
\draw [rotate around={-0.7073193685442621:(4.796742223721992,1.9905579037715853)},dotted,color=qqqqff] (4.796742223721992,1.9905579037715853) ellipse (0.7970491859738003cm and 0.445448289568427cm);
\draw [dash pattern=on 2pt off 2pt] (1.4,-0.4)-- (1.6,0.4);
\draw [dash pattern=on 2pt off 2pt] (1.4,-0.4)-- (1.8,0.4);
\draw [dotted,color=ffqqqq] (2.8,2.2)-- (2.8,1.8);
\draw [dotted,color=ffqqqq] (2.8,1.8)-- (3.8,1.8);
\draw [dotted,color=ffqqqq] (3.8,1.8)-- (3.8,2.2);
\draw [dotted,color=ffqqqq] (3.8,2.2)-- (2.8,2.2);
\draw [dotted] (5.4,2.)-- (5.6,2.2);
\draw (-0.4,1.2)-- (0.,1.6);
\draw (0.,1.6)-- (0.,2.);
\draw [dotted] (6.2,2.)-- (6.,1.8);
\draw [dotted] (5.4,2.)-- (5.6,1.8);
\draw [dotted] (5.6,1.8)-- (6.,1.8);
\draw [dotted] (5.6,2.2)-- (6.,2.2);
\draw [dotted] (6.,2.2)-- (6.2,2.);
\begin{scriptsize}
\draw [fill=qqqqff] (-1.2,2.) circle (2.5pt);
\draw [fill=qqqqff] (-1.2,1.6) circle (2.5pt);
\draw [fill=qqqqff] (-1.2,1.2) circle (2.5pt);
\draw [fill=qqqqff] (-1.2,0.8) circle (2.5pt);
\draw [fill=qqqqff] (-0.8,2.) circle (2.5pt);
\draw [fill=qqqqff] (-0.6,2.) circle (2.5pt);
\draw [fill=qqqqff] (-0.4,1.6) circle (2.5pt);
\draw [fill=qqqqff] (-0.4,2.) circle (2.5pt);
\draw [fill=qqqqff] (-0.2,2.) circle (2.5pt);
\draw [fill=qqqqff] (-0.4,1.2) circle (2.5pt);
\draw [fill=qqqqff] (0.,2.) circle (2.5pt);
\draw [fill=qqqqff] (-1.2,0.4) circle (2.5pt);
\draw [fill=qqqqff] (0.6,0.4) circle (2.5pt);
\draw [fill=qqqqff] (0.6,0.8) circle (2.5pt);
\draw [fill=qqqqff] (0.6,1.2) circle (2.5pt);
\draw [fill=qqqqff] (0.6,1.6) circle (2.5pt);
\draw [fill=qqqqff] (0.6,2.) circle (2.5pt);
\draw [fill=qqqqff] (1.,2.) circle (2.5pt);
\draw [fill=qqqqff] (1.4,1.6) circle (2.5pt);
\draw [fill=qqqqff] (1.2,2.) circle (2.5pt);
\draw [fill=qqqqff] (1.4,2.) circle (2.5pt);
\draw [fill=qqqqff] (1.6,2.) circle (2.5pt);
\draw [fill=qqqqff] (-0.6,-0.4) circle (2.5pt);
\draw [fill=qqqqff] (0.2,-0.4) circle (2.5pt);
\draw [fill=qqqqff] (-0.2,-0.8) circle (2.5pt);
\draw [fill=qqqqff] (1.4,-0.4) circle (2.5pt);
\draw [fill=qqqqff] (4.4,-0.8) circle (2.5pt);
\draw [fill=qqqqff] (4.4,-0.4) circle (2.5pt);
\draw [fill=qqqqff] (4.4,0.) circle (2.5pt);
\draw [fill=qqqqff] (4.4,0.4) circle (2.5pt);
\draw [fill=qqqqff] (3.4,0.8) circle (2.5pt);
\draw [fill=qqqqff] (4.2,0.8) circle (2.5pt);
\draw [fill=qqqqff] (5.,0.8) circle (2.5pt);
\draw [fill=qqqqff] (5.4,0.8) circle (2.5pt);
\draw [fill=qqqqff] (2.8,1.2) circle (2.5pt);
\draw [fill=qqqqff] (3.4,1.2) circle (2.5pt);
\draw [fill=qqqqff] (4.2,1.2) circle (2.5pt);
\draw [fill=qqqqff] (4.4,1.2) circle (2.5pt);
\draw [fill=qqqqff] (2.6,1.6) circle (2.5pt);
\draw [fill=qqqqff] (2.6,2.) circle (2.5pt);
\draw [fill=qqqqff] (3.,2.) circle (2.5pt);
\draw [fill=qqqqff] (3.2,2.) circle (2.5pt);
\draw [fill=qqqqff] (3.4,2.) circle (2.5pt);
\draw [fill=qqqqff] (4.2,1.6) circle (2.5pt);
\draw [fill=qqqqff] (4.2,2.) circle (2.5pt);
\draw [fill=qqqqff] (4.6,1.2) circle (2.5pt);
\draw [fill=qqqqff] (4.8,2.) circle (2.5pt);
\draw [fill=qqqqff] (5.,2.) circle (2.5pt);
\draw [fill=qqqqff] (4.6,2.) circle (2.5pt);
\draw [fill=qqqqff] (5.6,2.) circle (2.5pt);
\draw [fill=qqqqff] (5.8,2.) circle (2.5pt);
\draw [fill=qqqqff] (6.,2.) circle (2.5pt);
\draw [fill=qqqqff] (5.8,-0.4) circle (2.5pt);
\draw [fill=qqqqff] (5.2,-0.4) circle (2.5pt);
\draw [fill=qqqqff] (3.6,2.) circle (2.5pt);
\draw [fill=qqqqff] (0.,1.6) circle (2.5pt);
\end{scriptsize}
\end{tikzpicture}
\caption{In the first sum, $A_u$ and $A_v$ are independent under $P_\omega$. In the second sum, $\max_{1\le j \le k-2} M(u_j)$, $\max_{1\le j \le k-2} M(v_j)$, $M(u_{k-1}, v_{k-1})$ are represented by the rectangle, ellipse and hexagon respectively, and are independent under $P_\omega$ (for the definition of $M(u_{k-1}, v_{k-1})$, see \eqref{M(u,v)}). }
\label{f:fig2}
\end{figure}

\noindent We treat the two sums on the right-hand side. See Figure \ref{f:fig2}.

{\bf First sum:} $\sum_{k= \lambda_1 n+1}^n \E(\cdots)$. When $k>\lambda_1 n$, the events $A_u$ and $A_v$ are independent under $P_\omega$, so
\begin{eqnarray*}
    P_\omega(A_u \cap A_v)
 &=& P_\omega(A_u) \, P_\omega(A_v) 
    \\
 &\lesssim& \Big( p^2 n^{1+\alpha} \, (b+n)^\alpha \sum_{j=1}^{\lambda_1 n} m^{-(b+n)-j} \Lambda (u_j) \Big) \times
\\
&&
\Big( p^2 n^{1+\alpha} \, (b+n)^\alpha \sum_{j=1}^{\lambda_1 n} m^{-(b+n)-j} \Lambda (v_j) \Big) ,
\end{eqnarray*}

\noindent the inequality being a consequence of \eqref{P(Au)<}. We take the expectation with respect to $\P$ on both sides, to see that
\begin{eqnarray}
 &&\sum_{k=\lambda_1 n+1}^n \E \sum_{x\in \T_k} \; \sum_{(u, \, v)}  P_\omega(A_u \cap A_v)
    \nonumber
    \\
 &\lesssim& \sum_{k=\lambda_1 n+1}^n  m^{n-k} \; m^{2k} \,  \Big( p^2 n^{1+\alpha} \, (b+n)^\alpha \sum_{j=1}^{\lambda_1 n} m^{-(b+n)-j} m^j\Big)^2
    \nonumber
    \\ 
 &\asymp& \Big( \, p^2 n^{2+\alpha} (b+n)^\alpha m^{-b} \, \Big)^2
    \nonumber
    \\
 &\asymp & (\e(Z))^2, 
    \label{E(Z2)bis}
\end{eqnarray}

\noindent the last line being a consequence of \eqref{E(Z)bis}.

{\bf Second sum:} $\sum_{k=1}^{\lambda_1 n} \E(\cdots)$. This time, we argue differently, first by conditioning on $X(u)$ and $X(v)$. For $\ell_u$, $\ell_v \in [\lambda_1 n, \, \lambda_2n]$, we have, by \eqref{P(X_0=t)}, for $w=u$ or $v$,
$$
\p (X(w) = \ell_w)
=
c_{15} \, p(\ell_w +\gamma)^\alpha m^{-\ell_w}
\asymp
pn^\alpha m^{-\ell_w}.
$$

\noindent So, for $1\le k\le \lambda_1 n$,
\begin{eqnarray*}
    \E \sum_{x\in \T_k} \; \sum_{(u, \, v)} P_\omega(A_u \cap A_v)
 &\asymp& \E \sum_{x\in \T_k} \; \sum_{(u, \, v)} \sum_{\ell_u, \, \ell_v = \lambda_1 n}^{\lambda_2 n} p^2 n^{2\alpha} m^{-\ell_u-\ell_v}
    \\
 && P_\omega(A_u \cap A_v\, | \, X(u) = \ell_u, \, X(v) = \ell_v) .
\end{eqnarray*}

\noindent Write\footnote{Notation: $\max_\varnothing := 0$.}
\begin{eqnarray}
    B_{x,u,v}
 &:=& \max_{k\le j \le \lambda_1 n} M (x_j) \vee M(u_{k-1}, \, v_{k-1}),
    \label{B}
    \\
    M(u_{k-1}, \, v_{k-1})
 &:=& \max_{r \in \bro (u_{k-1}, \, v_{k-1})} \max_{w \in \T_0(r)} (X(w)-(k-1))^+\, ,
    \label{M(u,v)}
\end{eqnarray}

\noindent where $a\vee b := \max\{ a, \, b\}$, and $\bro (u_{k-1}, \, v_{k-1}) := \bro (u_{k-1}) \backslash \{ v_{k-1}\} = \bro (v_{k-1}) \backslash \{ u_{k-1}\}$. [Observe that $u_j =x_j =v_j$ for $k\le j \le \lambda_1 n$.] The random variables $X(u)$, $X(v)$, $\max_{1\le j\le k-2} M(u_j)$, $\max_{1\le j\le k-2} M(v_j)$ and $B_{x,u,v}$ are independent under $P_\omega$ (see Figure 2). As such, for $\ell_u$, $\ell_v \in [\lambda_1 n, \, \lambda_2n]$,
\begin{eqnarray*} 
 &&P_\omega(A_u \cap A_v\, | \, X(u) = \ell_u, \, X(v) = \ell_v)
    \\
 &=& P_\omega \Big( \max_{1\le j\le k-2} M(u_j) \vee B_{x,u,v} \ge b+n - \ell_u, \, 
    \\
 && \qquad\qquad \max_{1\le j\le k-2} M(v_j) \vee B_{x,u,v} \ge b+n - \ell_v \Big)
    \\
 &\le& P_\omega(B_{x,u,v} \ge  b+n - \ell_u) + P_\omega(B_{x,u,v} \ge  b+n - \ell_v) +
    \\
 && +
    P_\omega( \max_{1\le j\le k-2} M(u_j) \ge b+n - \ell_u) \, P_\omega( \max_{1\le j\le k-2} M(v_j)\ge b+n - \ell_v).
\end{eqnarray*}

\noindent The first two probability expressions on the right-hand side play the same role by symmetry in $\ell_u$ and $\ell_v$, so let us only look at the first one: By \eqref{M},
\begin{eqnarray*}
    P_\omega( B_{x,u,v} \ge  b+n - \ell_u) 
 &\lesssim&
p (b+n - \ell_u)^\alpha \sum_{j=k-1}^{\lambda_1 n} m^{-(b+n - \ell_u)-j} \Lambda (u_j)  
    \\
 &\lesssim&
p (b+n)^\alpha \sum_{j=k-1}^{\lambda_1 n} m^{-(b+n - \ell_u)-j} \Lambda (u_j) \, .
\end{eqnarray*}

\noindent [Note that $u_j =x_j$ for $k\le j\le \lambda_1 n$.] Similarly, for the third probability expression, we have, by \eqref{M} again, for $w=u$ or $v$, 
$$
P_\omega( \max_{1\le j\le k-2} M(w_j) \ge b+n - \ell_w)
\lesssim
p (b+n)^\alpha \sum_{j=1}^{k-2} m^{-(b+n - \ell_w)-j} \Lambda (w_j) \, .
$$

\noindent Assembling these pieces together yields, for $1\le k\le \lambda_1 n$,
\begin{eqnarray*}
 &&\E \sum_{x\in \T_k} \; \sum_{(u, \, v)} P_\omega(A_u \cap A_v)
    \\
 &\lesssim& \E \sum_{x\in \T_k} \; \sum_{(u, \, v)} \sum_{\ell_u, \, \ell_v= \lambda_1 n}^{\lambda_2 n} p^2 n^{2\alpha} m^{-\ell_u-\ell_v}\, p (b+n)^\alpha \sum_{j=k-1}^{\lambda_1 n} m^{-(b+n - \ell_u)-j} \Lambda (u_j)
    \\
 && + \E \sum_{x\in \T_k} \; \sum_{(u, \, v)} \sum_{\ell_u, \, \ell_v= \lambda_1 n}^{\lambda_2 n} p^2 n^{2\alpha} m^{-\ell_u-\ell_v} \, p^2 (b+n)^{2\alpha} 
    \\
 && \Big( \sum_{j=1}^{k-2} m^{-(b+n - \ell_u)-j} \Lambda(u_j) \Big) \Big( \sum_{j=1}^{k-2} m^{-(b+n - \ell_v)-j} \Lambda (v_j) \Big)\, .
\end{eqnarray*}

\noindent On the right-hand side, both expectations can be easily estimated by means of the branching property. The first expectation is
\begin{eqnarray*}
 &\lesssim& m^{n-k} \, m^{2k} \, \sum_{\ell_u, \, \ell_v= \lambda_1 n}^{\lambda_2 n} p^2 n^{2\alpha} m^{-(\ell_u+\ell_v)} \, p (b+n)^\alpha \sum_{j=k-1}^{\lambda_1 n} m^{-(b+n - \ell_u)-j} m^j
 \\
 &\asymp& m^{k-\lambda_1 n} \, p^3 n^{1+2\alpha} (b+n)^\alpha  m^{-b} (\lambda_1 n-k) \, ,
\end{eqnarray*}

\noindent whereas the second probability expression is 
\begin{eqnarray*}
 &\lesssim& m^{n-k} \, m^{2k} \, \sum_{\ell_u, \, \ell_v = \lambda_1 n}^{\lambda_2 n} p^2 n^{2\alpha} m^{-\ell_u -\ell_v} \, p^2 (b+n)^{2\alpha} 
    \\
 && \Big( \sum_{j=1}^{k-1} m^{-(b+n - \ell_u)-j}m^j \Big) \Big( \sum_{j=1}^{k-1} m^{-(b+n - \ell_v)-j}m^j \Big)
    \\
 &\asymp& k^2 m^{k-n-2b} p^4 n^{2+2\alpha} (b+n)^{2\alpha}.
\end{eqnarray*} 
 
\noindent Consequently, for $1\le k\le \lambda_1 n$,
\begin{eqnarray*}
 &&\E \sum_{x\in \T_k} \; \sum_{(u, \, v)} P_\omega(A_u \cap A_v)
    \\
 &\lesssim& m^{k-\lambda_1 n} \, p^3 n^{1+2\alpha} (b+n)^\alpha  m^{-b} (\lambda_1 n-k)
    +
    k^2 m^{k-n-2b} p^4 n^{2+2\alpha} (b+n)^{2\alpha} \, .
\end{eqnarray*}

\noindent Summing over $1\le k\le \lambda_1 n$ yields that
\begin{eqnarray*}
 &&\sum_{k=1}^{\lambda_1 n} \E \sum_{x\in \T_k} \; \sum_{(u, \, v)} P_\omega(A_u \cap A_v)
    \\
 &\lesssim& p^3 n^{1+2\alpha} (b+n)^\alpha  m^{-b}
    +
    m^{-(1-\lambda_1)n-2b} p^4 n^{3+2\alpha} (b+n)^{2\alpha}
    \\
 &\asymp& p^3 n^{1+2\alpha} (b+n)^\alpha  m^{-b} \, ,
\end{eqnarray*}
 
\noindent which, by \eqref{E(Z)bis}, is $\lesssim \e(Z)$ if $p n^{\alpha-1} \le 1$. Combining this with \eqref{E(Z2)bis} and \eqref{E(Z2)}, we see that $\e(Z^2) \lesssim \e(Z) + (\e Z)^2$, as long as $p n^{\alpha-1} \le 1$. 

Under the condition $p\, n^{1+\alpha}\le 1$, we have $\e (Z) \lesssim 1$ by \eqref{E(Z)bis}, so $(\e Z)^2 \lesssim \e(Z)$; it follows from the Cauchy--Schwarz inequality that $\p(Z\ge 1) \ge \frac{(\e Z)^2}{\e(Z^2)} \gtrsim \e (Z) \asymp p^2 \, n^{2+\alpha} (b+n)^\alpha \, m^{-b}$ (by \eqref{E(Z)bis}). The lemma follows now from \eqref{lb}.$\qed$

\subsection{Second step: The spinal advantage}

Let $\alpha>-2$ and $\varepsilon>0$. 

\medskip

\begin{lemma}
\label{l:E(Xn)>}

 Assume $\E(\nu^3) <\infty$ and $m:= \E (\nu)>1$.
 Under the assumption $\eqref{hyp}$, for any $K> \frac{2}{2+\alpha}$, there exist constants $s\ge K$ and $c_{19}>0$ such that for $0<p<1$ and $n\le (\frac1p)^{(s+1)/[(2+\alpha)s + \alpha]}$,
\begin{equation}
    \e(X_n) \ge c_{19} \, p^{s+1}\, n^{(2+\alpha)s + \alpha} .
    \label{E(Xn)>}
\end{equation}

\end{lemma}

\medskip

\noindent {\it Proof of Lemma \ref{l:E(Xn)>}.} [The condition $K> \frac{2}{2+\alpha}$ is to make sure that $(2+\alpha)s + \alpha>0$ whenever $s\ge K$.]

Let $i\ge 1$ be an integer. Let $a_i:= 2^i$ and $b_i:= (2+\alpha) (2^i -1)$ (which explains the condition $\alpha>-2$: so that $b_i >0$). Applying Lemma \ref{l:iteration} $i$ times, we see that for any integer $i\ge 1$ and any constant $c>0$, there exists a constant $c(i) >0$ such that for $0<p<1$ and $n\ge 1$ with $p^{a_i} n^{b_i+\alpha} \le c$, we have
\begin{equation}
    \p(X_{in} \ge b) \ge c(i) \,  p^{a_i} n^{b_i}  (b+n)^\alpha \,  m^{-b}, 
    \qquad \forall b>0.
    \label{induction_i}
\end{equation}

\noindent [For $i=1$, \eqref{induction_i} follows from Lemma \ref{l:iteration} with $\gamma=1$ and $t=b$. Assuming \eqref{induction_i} holds for $i$, it is immediately seen to hold for $i+1$: It suffices to apply Lemma \ref{l:iteration} to $\gamma=n$, $t=b$ and to $p^{a_i}n^{b_i}$ in place of $p$.] 

For integers $\ell \in [0, \, n]$, we use the above inequality for $\p(X_{(i-1)n} \ge b)$, and apply Lemma \ref{l:iteration} to $n+\ell$ in place of $n$, to see that there exists a constant $c'(i) >0$ for $0<p<1$ and $n\ge 1$ with $p^{a_i} (n+\ell)^{b_i+\alpha} \le c$, and for all integers $\ell \in [0, \, n]$, 
$$
\p(X_{in+\ell} \ge b) \ge c'(i) \,  p^{a_i} n^{b_i}  (b+n)^\alpha \,  m^{-b}, 
\qquad \forall b>0.
$$

\noindent Integrating over $b$, this yields the existence of a constant $c''(i) >0$, depending on $i$, such that for $0<p<1$ and $n\ge 1$ with $p^{a_i} n^{b_i+\alpha} \le c$,
$$
\e(X_{in+\ell}) \ge c''(i) \, p^{a_i} n^{b_i+\alpha} .
$$

\noindent We choose (and fix) $i$ sufficiently large, how large depending on $\alpha$, such that $a_i \ge K+1$. The lemma follows with $s:= a_i -1$.\qed

\bigskip 

The rest of the proof of the lower bound in Theorem \ref{t:main} consists in improving the lower bound for $\e(X_n)$ in \eqref{E(Xn)>}, and making it (strictly) greater than $\frac{1}{m-1}$, so that by virtue of the first inequality in \eqref{encadrement_F}, which says that $F_\infty \ge \frac{\e(X_n)-\frac{1}{m-1}}{m^n}$, it will give the desired lower bound for the free energy $F_\infty$ as stated in Theorem \ref{t:main}.

Without loss of generality, we assume that the law of $X_0$, conditionally on $X_0>0$, is absolutely continuous.

To improve \eqref{E(Xn)>}, we start with a new lower bound for $X_n$. For any vertex $v\in \T^{(n)}$, we write $X(v)$ for the random variable associated with the vertex $v$: so if $|v|=j$, then $X(v)$ is distributed as $X_j$. Let $1\le k <\ell < n$ be integers; the values of $k$ and $\ell$, both depending on $(n, \, p)$, will be given later. For $u\in \T_0$, let
$$
M^*_k(u)
:=
\max_{v\in \T_0(u_k) \setminus \{ u\}} X(v),
\qquad
N^*_k(u)
:=
\max_{v\in \T_0 \setminus \T_0(u_k)} X(v) .
$$

\noindent [So $M^*_k(u) \vee N^*_k(u) \vee X(u)$ coincides with $\max_{w\in \T_0} X(w)$; moreover, $M^*_k(u)$, $N^*_k(u)$, $X(u)$ are independent random variables under $P_\omega$.] 

Assume there exists $u\in \T_0$ such that $X(u) > N^*_k(u)$ and $M^*_k(u) \le \ell<X(u)$. If such a vertex $u$ exists (which must be unique, by definition), $X_n \ge X(u) - n + \sum_{j=0}^\ell \sum_{v\in \bro(u_j)} X(v)$, which is greater than $\ell - n + \sum_{j=0}^{k-1} \sum_{v\in \bro(u_j)} X(v)$.\footnote{Strictly speaking, we should write $X(\mathfrak{e}_n)$ in place of $X_n$.} [In case $\ell - n + \sum_{j=0}^{k-1} \sum_{v\in \bro(u_j)} X(v)$ is negative, the statement is, of course, trivial.] We arrive at the following inequality:
\begin{eqnarray}
    X_n
 &\ge& \sum_{u\in \T_0} 
    {\bf 1}_{\{ X(u) > N^*_k(u) \vee \ell \} } \, 
    {\bf 1}_{\{ M^*_k(u) \le \ell\} } \,
    \Big( \ell - n + \sum_{j=0}^{k-1} \sum_{v\in \bro(u_j)} X(v) \Big) .
    \nonumber
    \\
 &=:& \sum_{u\in \T_0} 
    {\bf 1}_{\{ X(u) > N^*_k(u) \vee \ell \} } \, \xi(u),
    \label{X>}
\end{eqnarray}

\noindent where, for $u\in \T_0$,
$$
\xi (u)
:=
{\bf 1}_{\{ M^*_k(u) \le \ell\} } \,
\Big( \ell - n + \sum_{j=0}^{k-1} \sum_{v\in \bro(u_j)} X(v) \Big) .
$$

\begin{figure}[h]

\definecolor{ttttff}{rgb}{0.2,0.2,1.}
\definecolor{qqqqff}{rgb}{0.,0.,1.}
\begin{tikzpicture}[line cap=round,line join=round,>=triangle 45, scale=1, every node/.style={transform shape}, x=4.0cm,y=3.5cm]
\clip(-1.6482839233060453,0.2260138959450558) rectangle (2.2035203804636376,2.751869859542719);
\draw (1.,2.)-- (1.,1.8);
\draw [line width=2.8pt] (0.2,1.8)-- (0.,2.);
\draw [line width=2.8pt] (0.4,1.4)-- (0.4,1.2);
\draw [line width=2.8pt] (0.4,1.2)-- (0.4,1.);
\draw [line width=2.8pt] (0.4,1.)-- (0.4,0.8);
\draw [line width=2.8pt] (0.4,0.8)-- (0.4,0.6);
\draw (0.4,1.4)-- (1.,1.6);
\draw (1.,1.6)-- (1.,1.8);
\draw [line width=2.8pt] (0.4,1.4)-- (0.4,1.6);
\draw [line width=2.8pt] (0.4,1.6)-- (0.2,1.8);
\draw (0.41878383133475316,0.6444690755430719) node[anchor=north west] {$\mathfrak{e}_n$};
\draw (0.41878383133475316,1.441046405380259) node[anchor=north west] {$\mathfrak{e}_k$};
\draw (0.41878383133475316,1.24442288725589) node[anchor=north west] {$\mathfrak{e}_{k+1}$};
\draw (0.41878383133475316,0.8461342223372964) node[anchor=north west] {$\mathfrak{e}_{n-1}$};
\draw (0.41878383133475316,1.647753180844339) node[anchor=north west] {$\mathfrak{e}_{k-1}$};
\draw (-0.1307536936795079,2.041000217093077) node[anchor=north west] {$\mathfrak{e}_2$};
\draw (-1.,2.2)-- (-1.,2.4);
\draw (-0.8,2.2)-- (-0.8,2.4);
\draw (-0.6,2.2)-- (-0.6,2.4);
\draw (1.,2.)-- (1.,2.2);
\draw (1.,2.2)-- (1.,2.4);
\draw (0.4,2.2)-- (0.6,2.4);
\draw [line width=2.8pt] (0.,2.)-- (0.,2.2);
\draw [line width=2.8pt] (0.,2.2)-- (0.,2.4);
\draw (-0.6,2.2)-- (-0.6,2.4);
\draw (-0.1307536936795079,2.2477069925571573) node[anchor=north west] {$\mathfrak{e}_1$};
\draw (-0.2,0.8)-- (0.4,0.6);
\draw (0.8,2.)-- (0.8,2.2);
\draw (0.8,2.2)-- (0.8,2.4);
\draw (0.4,1.6)-- (0.6,1.8);
\draw (0.6,1.8)-- (0.8,2.);
\draw (0.2,1.8)-- (0.4,2.);
\draw (0.4,2.)-- (0.4,2.2);
\draw (0.4,0.6)-- (-0.8,0.8);
\draw [dash pattern=on 4pt off 4pt] (-0.8,0.8)-- (-1.,2.2);
\draw [dash pattern=on 4pt off 4pt] (-0.2,0.8)-- (-0.8,2.2);
\draw (0.6,2.59) node[anchor=north west] {\scriptsize$E_\omega(\xi(\mathfrak{e}_0))$};
\draw (-0.9827889388851053,2.59) node[anchor=north west] {\scriptsize$P_\omega\{ X(\mathfrak{e}_0) > N^*_k(\mathfrak{e}_0) \vee \ell \}$};
\draw [dash pattern=on 4pt off 4pt] (0.4,1.2)-- (-0.6,2.2);
\draw (-0.6,2.2)-- (-0.4,2.4);
\draw (1.,2.2)-- (1.2,2.4);
\draw (0.,2.2)-- (0.4,2.4);
\draw (-0.1307536936795079,2.45) node[anchor=north west] {$\mathfrak{e}_0$};
\draw [rotate around={-0.8107352192623714:(0.7782519554954525,2.4201306930662194)},dotted,color=ttttff] (0.7782519554954525,2.4201306930662194) ellipse (2.2080723699913185cm and 0.48332852271109616cm);
\draw [dotted] (-1.1556336677052708,2.6)-- (0.11126733541054135,2.6);
\draw [dotted] (0.11126733541054135,2.6)-- (0.11126733541054135,2.288732664589459);
\draw [dotted] (0.11126733541054135,2.288732664589459)-- (-1.1556336677052708,2.288732664589459);
\draw [dotted] (-1.1556336677052708,2.288732664589459)-- (-1.1556336677052708,2.6);
\begin{scriptsize}
\draw [fill=qqqqff] (-1.,2.2) circle (2.5pt);
\draw [fill=qqqqff] (-0.8,2.2) circle (2.5pt);
\draw [fill=qqqqff] (1.,2.) circle (2.5pt);
\draw [fill=qqqqff] (0.,2.) circle (2.5pt);
\draw [fill=qqqqff] (1.,1.8) circle (2.5pt);
\draw [fill=qqqqff] (0.2,1.8) circle (2.5pt);
\draw [fill=qqqqff] (0.4,0.6) circle (2.5pt);
\draw [fill=qqqqff] (0.4,1.4) circle (2.5pt);
\draw [fill=qqqqff] (0.4,1.2) circle (2.5pt);
\draw [fill=qqqqff] (0.4,1.) circle (2.5pt);
\draw [fill=qqqqff] (0.4,0.8) circle (2.5pt);
\draw [fill=qqqqff] (1.,1.6) circle (2.5pt);
\draw [fill=qqqqff] (0.4,1.6) circle (2.5pt);
\draw [fill=qqqqff] (-1.,2.2) circle (2.5pt);
\draw [fill=qqqqff] (-0.8,2.2) circle (2.5pt);
\draw [fill=qqqqff] (-1.,2.4) circle (2.5pt);
\draw [fill=qqqqff] (-0.8,2.4) circle (2.5pt);
\draw [fill=qqqqff] (-0.6,2.4) circle (2.5pt);
\draw [fill=qqqqff] (-0.6,2.2) circle (2.5pt);
\draw [fill=qqqqff] (-0.6,2.4) circle (2.5pt);
\draw [fill=qqqqff] (1.,2.2) circle (2.5pt);
\draw [fill=qqqqff] (1.,2.4) circle (2.5pt);
\draw [fill=qqqqff] (0.4,2.2) circle (2.5pt);
\draw [fill=qqqqff] (0.6,2.4) circle (2.5pt);
\draw [fill=qqqqff] (0.4,2.4) circle (2.5pt);
\draw [fill=qqqqff] (0.,2.2) circle (2.5pt);
\draw [fill=qqqqff] (0.,2.4) circle (2.5pt);
\draw [fill=qqqqff] (0.8,2.4) circle (2.5pt);
\draw [fill=qqqqff] (-0.2,0.8) circle (2.5pt);
\draw [fill=qqqqff] (0.8,2.) circle (2.5pt);
\draw [fill=qqqqff] (0.8,2.2) circle (2.5pt);
\draw [fill=qqqqff] (0.6,1.8) circle (2.5pt);
\draw [fill=qqqqff] (0.4,2.) circle (2.5pt);
\draw [fill=qqqqff] (-0.8,0.8) circle (2.5pt);
\draw [fill=qqqqff] (-0.6,2.2) circle (2.5pt);
\draw [fill=qqqqff] (-0.4,2.4) circle (2.5pt);
\draw [fill=qqqqff] (1.2,2.4) circle (2.5pt);
\end{scriptsize}
\end{tikzpicture}

\caption{The random variables $P_\omega\{ X(\mathfrak{e}_0) > N^*_k(\mathfrak{e}_0) \vee \ell \}$ and $E_\omega (\xi (\mathfrak{e}_0))$ are independent under $\Q$.}
\label{f:fig3}

\end{figure}

\noindent Note that $X(u)$, $N^*_k(u)$ and $\xi(u)$ are independent under $P_\omega$. Hence
$$
    E_\omega(X_n)
    \ge
    \sum_{u\in \T_0} 
    P_\omega\{ X(u) > N^*_k(u) \vee \ell \} \, 
    E_\omega (\xi (u))\, ,
$$

\noindent Taking expectation with respect to $\P$, we obtain, by the many-to-one formula \eqref{many-to-one},
$$
    \e(X_n)
    \ge
    m^n \, \E_\Q \Big[ P_\omega\{ X(\mathfrak{e}_0) > N^*_k(\mathfrak{e}_0) \vee \ell \} \, 
    E_\omega (\xi (\mathfrak{e}_0)) \Big] \, .
$$

\noindent By the spinal decomposition theorem, the random variables $P_\omega\{ X(\mathfrak{e}_0) > N^*_k(\mathfrak{e}_0) \vee \ell \}$ and $E_\omega (\xi (\mathfrak{e}_0))$ are independent under $\Q$. See Figure \ref{f:fig3}. So
\begin{equation}
    \e(X_n)
    \ge
    m^n \, (\Q \otimes P_\omega) \{ X(\mathfrak{e}_0) > N^*_k(\mathfrak{e}_0) \vee \ell \} \; 
    (\E_\Q \otimes E_\omega) (\xi (\mathfrak{e}_0)) \, .
    \label{xi}
\end{equation}

We study $(\E_\Q \otimes E_\omega) (\xi (\mathfrak{e}_0))$ on the right-hand side. Since $\ell <n$, we have,
\begin{eqnarray*}
    \xi (\mathfrak{e}_0)
 &\ge& \ell - n + {\bf 1}_{\{ M^*_k(\mathfrak{e}_0) \le \ell\} } \, \sum_{j=0}^{k-1} \sum_{v\in \bro(\mathfrak{e}_j)} X(v) 
    \\
 &=& \ell - n + \sum_{j=0}^{k-1} \sum_{v\in \bro(\mathfrak{e}_j)} X(v) 
    -
    \sum_{j=0}^{k-1} {\bf 1}_{\{ M^*_k(\mathfrak{e}_0) > \ell\} } \sum_{v\in \bro(\mathfrak{e}_j)} X(v) \, .
\end{eqnarray*}

\noindent We take expectation with respect to $(\Q\otimes P_\omega)$ on both sides. By the spinal decomposition theorem,
\begin{eqnarray}
    (\E_\Q \otimes E_\omega) (\xi(\mathfrak{e}_0))
 &\ge& \ell - n + \sum_{j=0}^{k-1} \e(X_j) \, \E_\Q[\# \bro(\mathfrak{e}_0) ]
    \nonumber
    \\ 
 && -
    \sum_{j=0}^{k-1} (\E_\Q \otimes E_\omega ) \Big( {\bf 1}_{\{ M^*_k(\mathfrak{e}_0) > \ell\} } \sum_{v\in \bro(\mathfrak{e}_j)} X(v) \Big) \, .
    \label{E(xi)}
\end{eqnarray}

Let us have a closer look at the last $(\E_\Q \otimes E_\omega ) (\cdots )$ expression on the right-hand side. By the trivial inequality $X(v) \le \sum_{r\in \T_0(v)} X(r)$, we have 
$$
\sum_{v\in \bro(\mathfrak{e}_j)} X(v) 
\le 
\sum_{v\in \bro(\mathfrak{e}_j)} \sum_{r\in \T_0(v)} X(r) 
\le  
\sum_{r\in \T_0(\mathfrak{e}_{j+1})} X(r) \, .
$$

\noindent So by the Cauchy--Schwarz inequality,
\begin{eqnarray*}
 &&(\E_\Q \otimes E_\omega ) \Big( {\bf 1}_{\{ M^*_k(\mathfrak{e}_0) > \ell\} } \sum_{v\in \bro(\mathfrak{e}_j)} X(v) \Big)
    \\
 &\le&
[ (\Q \otimes P_\omega )\{ M^*_k(\mathfrak{e}_0) > \ell\}]^{1/2} \,
\Big\{ (\E_\Q \otimes E_\omega ) \Big[ \Big( \sum_{r\in \T_0(\mathfrak{e}_{j+1})} X(r) \Big)^{\! 2}\, \Big] \Big\}^{1/2} .
\end{eqnarray*}

\noindent By definition, $P_\omega \{ M^*_k(\mathfrak{e}_0) > \ell\} \le (\# \T_0(\mathfrak{e}_k)) \, \p (X_0 > \ell)$, which, by the assumption \eqref{hyp}, is $\lesssim (\# \T_0(\mathfrak{e}_k)) \, p \, \ell^\alpha \, m^{-\ell}$. Hence 
$$
(\Q \otimes P_\omega )\{ M^*_k(\mathfrak{e}_0) > \ell\} 
\le 
\E_\Q (\# \T_0(\mathfrak{e}_k)) \, \p (X_0 > \ell)
\lesssim
\E_\Q (\# \T_0(\mathfrak{e}_k)) \, p \, \ell^\alpha \, m^{-\ell}\, .
$$

\noindent By \eqref{E(T(ek))}, this yields
$$
(\Q \otimes P_\omega )\{ M^*_k(\mathfrak{e}_0) > \ell\} 
\lesssim
m^k \, p \, \ell^\alpha \, m^{-\ell}\, .
$$

\noindent On the other hand, 
$$
(\E_\Q \otimes E_\omega ) \Big[ \Big( \sum_{r\in \T_0(\mathfrak{e}_{j+1})} X(r) \Big)^{\! 2}\, \Big]
=
\E_\Q \Big[ \Big( \sum_{r\in \T_0(\mathfrak{e}_{j+1})} \e(X_0) \Big)^{\! 2}\, \Big]
+
\E_\Q \Big[ \sum_{r\in \T_0(\mathfrak{e}_{j+1})} \sigma^2 \Big] ,
$$

\noindent where $\sigma^2 := \mathrm{Var} (X_0) \le c_{20} \, p$, and $\e(X_0) = c_{21} \, p$. As such,
\begin{eqnarray*}
 &&(\E_\Q \otimes E_\omega ) \Big[ \Big( \sum_{r\in \T_0(\mathfrak{e}_{j+1})} X(r) \Big)^{\! 2}\, \Big]
    \\
 &\le& (c_{21} \, p)^2 \, \E_\Q[(\# \T_0(\mathfrak{e}_{j+1}))^2]
    +
    c_{20} \, p\, \E_\Q (\# \T_0(\mathfrak{e}_{j+1}) 
    \lesssim
    p^2 \, m^{2j} + p \, m^j \, ,
\end{eqnarray*}

\noindent the last inequality following from \eqref{E(T(ek))} and \eqref{E(T(ek)2)}. Consequently,
$$
(\E_\Q \otimes E_\omega ) \Big( {\bf 1}_{\{ M^*_k(\mathfrak{e}_0) > \ell\} } \sum_{v\in \bro(\mathfrak{e}_j)} X(v) \Big)
\lesssim
\Big[ m^k \, p \, \ell^\alpha \, m^{-\ell} (p^2 \, m^{2j} + p \, m^j) \Big]^{1/2} .
$$

\noindent As such, as long as we take
\begin{equation}
    k := \lfloor \frac{\ell}{4} \rfloor,
    \qquad
    \ell := n - c_{22} \, \log n \, ,
    \label{ell}
\end{equation}

\noindent we have, for some $c_{23} >0$,
$$
\sum_{j=0}^{k-1}
(\E_\Q \otimes E_\omega ) \Big( {\bf 1}_{\{ M^*_k(\mathfrak{e}_0) > \ell\} } \sum_{v\in \bro(\mathfrak{e}_j)} X(v) \Big)
\le
c_{23}.
$$

\noindent Going back to \eqref{E(xi)}, we obtain, with $c_{24} := \E_\Q[\# \bro(\mathfrak{e}_0) ] <\infty$,
$$
(\E_\Q \otimes E_\omega) (\xi(\mathfrak{e}_0))
\ge
\ell - n + c_{24} \sum_{j=0}^{k-1} \e(X_j) - c_{23}\, .
$$

Let $\varepsilon>0$. Let $s\ge K$ be the constants in Lemma \ref{l:E(Xn)>}. We choose $K$ so large that $\frac{s+1}{(2+\alpha)s + \alpha} \le \frac{1}{2+\alpha}+\varepsilon$. If we take 
\begin{equation}
    n = n(p)
    := 
    \Big( \frac1p \Big)^{\! \frac{s+1}{(2+\alpha)s + \alpha}} 
    \le 
    \Big( \frac1p \Big)^{\! \frac{1}{2+\alpha}+\varepsilon},
    \label{n}
\end{equation}

\noindent (the last inequality holding for all sufficiently small $p$), then by Lemma \ref{l:E(Xn)>}, $\sum_{j=0}^{\lfloor \frac{\ell}{4} \rfloor -1} \e(X_j) \ge c_{25} \, n$ for some constant $c_{25}>0$. With our choice of $\ell$ in \eqref{ell}, this implies $(\E_\Q \otimes E_\omega)(\xi(\mathfrak{e}_0)) \ge c_{26} \, n$, with $c_{26}>0$. Going back to \eqref{xi}, we obtain, for all sufficiently small $p>0$,
\begin{equation}
    \e(X_n)
    \ge
    m^n \, (\Q \otimes P_\omega) \{ X(\mathfrak{e}_0) > N^*_k(\mathfrak{e}_0) \vee \ell \} \, c_{25} \, n\, .
    \label{E(Xn)>bis}
\end{equation}

\noindent By the many-to-one formula \eqref{many-to-one} again, this yields
\begin{eqnarray*}
    \e(X_n)
 &\ge& \E \Big( \sum_{u\in \T_0} P_\omega \{ X(u) > N^*_k(u) \vee \ell \} \Big) \, c_{25} \, n
    \\
 &\ge& \E \Big( \sum_{u\in \T_0} {\bf 1}_{\{ c_{27} \, m^n \le \# \T_0 \le c_{28} \, m^n\} } \, P_\omega \{ X(u) > N^*_k(u) \vee \ell \} \Big) \, c_{25} \, n \, ,
\end{eqnarray*}

\noindent for any constants $c_{28}>c_{27}>0$. For any $u\in \T_0$, since the conditional law of $X_0$ given $X_0>0$ is assumed to be absolutely continuous, we have 
$$
P_\omega \{ X(u) > N^*_k(u) \vee \ell \} 
\ge 
P_\omega \{ X(u) = \max_{r\in \T_0} X(r) \} 
-
P_\omega \{ \max_{r\in \T_0} X(r) \le \ell\} .
$$

\noindent We have $\sum_{u\in \T_0} P_\omega \{ X(u) = \max_{r\in \T_0} X(r) \} = 1$, whereas $P_\omega \{ \max_{r\in \T_0} X(r) \le \ell\} = [1-\p(X_0 > \ell)]^{\# \T_0} \le [1-c_3 \ell^\alpha m^{-\ell}]^{\# \T_0}$ by \eqref{hyp}; hence
\begin{eqnarray*}
 &&\E \Big( \sum_{u\in \T_0} {\bf 1}_{\{ c_{27} \, m^n \le \# \T_0 \le c_{28} \, m^n\} } \, P_\omega \{ X(u) > N^*_k(u) \vee \ell \} \Big)
    \\
 &\ge& \P(c_{27} \, m^n \le \# \T_0 \le c_{28} \, m^n) - \E \Big( \sum_{u\in \T_0} {\bf 1}_{\{ c_{27} \, m^n \le \# \T_0 \le c_{28} \, m^n\} } [1-c_3 \ell^\alpha m^{-\ell}]^{\# \T_0} \Big)
    \\
 &=&\P(c_{27} \, m^n \le \# \T_0 \le c_{28} \, m^n) - \E \Big( (\# \T_0)\, {\bf 1}_{\{ c_{27} \, m^n \le \# \T_0 \le c_{28} \, m^n\} } [1-c_3 \ell^\alpha m^{-\ell}]^{\# \T_0} \Big).
\end{eqnarray*}

\noindent On the event $\{ c_{27} \, m^n \le \# \T_0 \le c_{28} \, m^n\}$, we have $(\# \T_0)\, [1-c_3 \ell^\alpha m^{-\ell}]^{\# \T_0} \le c_{28} \, m^n \exp( - c_3 \ell^\alpha m^{-\ell}\, c_{27} \, m^n)$, which tends to $0$ with the choice of $\ell := n- c_{22} \, \log n$, as long as $c_{22}>\frac{1-\alpha}{\log m}$. On the other hand, the constants $c_{28}>c_{27}>0$ can be chosen such that $\P(c_{27} \, m^n \le \# \T_0 \le c_{28} \, m^n) \to c_{29}>0$, $n\to \infty$. Hence
$$
\liminf_{n\to \infty} \E \Big( \sum_{u\in \T_0} {\bf 1}_{\{ c_{27} \, m^n \le \# \T_0 \le c_{28} \, m^n\} } \, P_\omega \{ X(u) > N^*_k(u) \vee \ell \} \Big)
\ge
c_{29}>0 \, .
$$

\noindent In view of \eqref{E(Xn)>bis}, we obtain: for some constant $c_{30}>0$ and all sufficiently small $p$, with $n=n(p)$ given in \eqref{n},
$$
\e(X_n)
\ge
c_{30} \, n .
$$

\noindent By the first inequality in \eqref{encadrement_F}, we get $F_\infty \ge \frac{\e(X_n)-\frac{1}{m-1}}{m^n} \ge \frac{c_{30} \, n -\frac{1}{m-1}}{m^n}$. The definition of $n$ in \eqref{n} yields that for an arbitrary $\varepsilon>0$ and all sufficiently small $p$, $F_\infty \ge \exp( - (\frac1p)^{\frac{1}{2+\alpha} + \varepsilon})$, proving the lower bound in Theorem \ref{t:main}.

[We mention that the lower bound is proved under the assumption $\E( \nu^3) <\infty$, instead of $\E(t^\nu)<\infty$ for some $t>1$.]\qed

\section{Proof of Theorem \ref{t:power_law}: lower bound}
\label{s:power_law}

We use the obvious stochastic inequality that $X_n$ is stochastically greater than or equal to $\max_{u\in \T_0} X(u) - n$. Hence for all $b\ge 0$,
\begin{eqnarray*}
    P_\omega (X_n >b)
 &\ge& P_\omega \Big( \max_{u\in \T_0} X(u) - n >b\Big)
    \\
 &=& 1- [1- \p(X_0 \le n+b)]^{\# \T_0} \, .
\end{eqnarray*}

\noindent By assumption, for all sufficiently large $n$ (say $n\ge n_1$), $\p (X_0 \le n+b) \ge p \, c_1 \, \ee^{-\theta (n+b)}$. Thus, for some constant $c_{31}>0$, all $n\ge n_1$ and all $b\ge 0$,
\begin{eqnarray*}
    P_\omega (X_n >b)
 &\ge& 1- \Big( 1- p\, c_1 \, \ee^{-\theta (n+b)}\Big)^{\# \T_0}
    \\
 &\ge& c_{31} \, \min\{ p\, \ee^{-\theta (n+b)}\# \T_0, \, 1\} \, .
\end{eqnarray*}

We take $n = N(p) := \lceil \frac{c_{32}+\log (1/p)}{(\log m)-\theta} \rceil$, where $c_{32}>0$ is a sufficiently large constant. Note that $p\, \ee^{-\theta (N(p)+b)} m^{N(p)} \ge 1$ if (and only if) $b \le B(p)$ where $B(p) := \frac{1}{\theta} [((\log m)-\theta)N(p) - \log (\frac1p)+\theta]$. Therefore, for all sufficiently small $p$ (such that $N(p)\ge n_1$) and all $0<b \le B(p)$, we have\footnote{To ensure that $\P(\# \T_0 \ge m^n)$ is greater than a positive constant, uniformly in $n$, it suffices to have $\E( \nu \log^+\! \nu) := \sum_{k=1}^\infty k(\log k) \P(\nu =k) <\infty$; see \cite{kesten-stigum} or \cite{lyons-pemantle-peres}.}
$$
\p(X_{N(p)} >b)
\ge
c_{31} \, \P(\# \T_0 \ge m^{N(p)})
=
c_{31} \, \P(\# \T_0 \ge m^n)
\ge
c_{33}\, ,
$$

\noindent where $c_{33}>0$ does not depend on $p$. This implies that $\E(X_n) = \int_0^\infty \P(X_n >b) \d b \ge c_{33}\, B(p)$, which, by definition of $B(p)$ and $N(p)$, is greater than $\frac{c_{33}\, (c_{32}+\theta)}{\theta}$. The latter is greater than $2$ if we choose $c_{32} := \frac{2\theta}{c_{33}}$. By the first inequality in \eqref{encadrement_F}, we get, for all sufficiently small $p$,
$$
F_\infty(p)
\ge
\frac{2-\frac{1}{m-1}}{m^{N(p)}}
=
(2-\frac{1}{m-1}) \exp \Big( - \lceil \frac{c_{32}+\log (1/p)}{(\log m)-\theta} \rceil \log m\Big),
$$

\noindent proving the lower bound in Theorem \ref{t:power_law}.

[The lower bound only requires $\E( \nu \log^+\! \nu) <\infty$, instead of $\E(t^\nu)<\infty$ for some $t>1$.]\qed

\section{Comments and questions}
\label{s:final}

We present some remarks and open problems.

\bigskip

{\bf (a) Change of measures?} Theorem A and Proposition \ref{p:p_c>0} in the introduction reveal the importance of $\E(X_0 \, m^{X_0})$. The latter strongly indicates that there could be a change-of-measures story hidden in the model. 

\medskip

\begin{problem}
\label{pr:change_of_measures}

 Is it possible to prove Theorem A by means of a change-of-measures argument? 
 
\end{problem}

\bigskip

{\bf (b) About the value of $p_c$ for non integer valued distributions.} Theorem A in the introduction gives the value of $p_c$ when $\nu = 2$ a.s.\ and $Y_0$ is integer valued. [They are valid whenever $\nu$ is deterministic (i.e., $\nu =m$ a.s.), with $\E[(Y_0-1)2^{Y_0}]$ and $\E(Y_0\, 2^{Y_0})$ replaced by $\E[((m-1)Y_0-1)m^{Y_0}]$ and $\E(Y_0\, m^{Y_0})$, respectively.] The following problem looks important to us.

\medskip

\begin{problem}
\label{pr:pc=?}

 Assume $\nu$ is deterministic. What can be said about $p_c$ without the assumption that $Y_0$ is integer valued?
 
\end{problem}

\medskip

Problem \ref{pr:pc=?}, which is borrowed from Derrida and Retaux~\cite{derrida-retaux}, seems challenging. For example, even assuming that $\nu = 2$ a.s., and that $Y_0$ takes values in $\{ \frac12, \, 1, \, \frac32, \, 2, \cdots\}$, we do not know what the value of $p_c$ should be in general.

\bigskip

{\bf (c) More about the value of $p_c$.} When $\nu$ is not deterministic, even if assuming $Y_0$ is integer valued as in Theorem A (see the introduction), it is not clear what $p_c$ should be. It is possible to have some bounds, but it seems to be hard to have an analytical expression.

\medskip

\begin{problem}
\label{prob:pc_GW}

 Assume $\nu$ is not deterministic and $Y_0$ takes values in $\{ 1, \, 2, \, \cdots\}$. What can be said about $p_c$?
 
\end{problem}

\bigskip

\noindent {\Large\bf Acknowledgements} 

\medskip

We are grateful to Bernard Derrida who introduced us to the problem, and with whom we have had regular discussions for two years. We wish to thank Nina Gantert for many discussions, Quentin Berger for enlightenment on renormalisation models, and Chunhua Ma, Bastien Mallein and Quan Shi for pointing out \cite{goldschmidt-przykucki} to us. Two anonymous referees have carefully read the manuscript; their insightful comments have led to improvements in the paper. The project was partly supported by ANR MALIN (ANR-16-CE93-0003); Y.H.\ also acknowledges support  from ANR SWiWS (ANR-17-CE40-0032).

\end{document}